# SEPARATION OF SEMIALGEBRAIC SETS

F. ACQUISTAPACE, C. ANDRADAS, AND F. BROGLIA

ABSTRACT. In this paper we study the problem of deciding whether two disjoint semialgebraic sets of an algebraic variety over $\mathbb{R}$ are separable by a polynomial. For that we isolate a dense subfamily of Spaces of Orderings, named Geometric, which suffice to test separation and that reduce the problem to the study of the behaviour of the semialgebraic sets in their boundary. Then we derive several characterizations for the generic separation, among which there is a Geometric Criterion that can be tested algorithmically. Finally we show how to check recursively whether we can pass from the generic separation to the separation of the two sets, yielding a decision procedure to solve the problem.

## INTRODUCTION

Let $M \subset \mathbb{R}^n$ be an algebraic variety, and $A, B \subset M$ two disjoint semialgebraic sets. When does it exist a regular function $f$ on $M$ which separates $A$ and $B$, i.e. such that $f(A) > 0$ and $f(B) < 0$?

We will show that this question is decidable, in the sense that there exists a decision procedure (theoretical algorithm) which takes $M, A$ and $B$ as input and produces YES or NO as output, according to whether $A$ and $B$ can be separated.

It is very easy to find examples of semialgebraic sets which cannot be separated by polynomials. The following two are among the more simple ones:

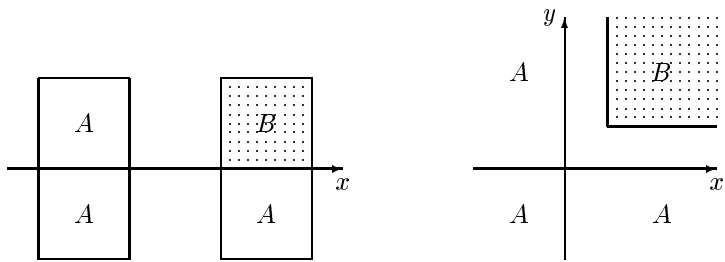

Although they seem different they are in fact quite similar. Indeed, consider, for the second, the one point compactification of $\mathbb{R}^2$. The shapes of $A$ and $B$ in a neighbourhood of the point at infinite look like the left hand side of the picture below. Now, if we blow-up this point we get the configuration of the right hand side, which obeys the same pattern that the first example above: there is an irreducible component of the boundary (the $x$–axis in the first example, the $y$–axis in the second) which has the property that any possible function separating $A$ and $B$

This work is partially supported by the EC contract CHRX-CT94-0506.
First and third authors are members of GNSAGA of CNR, and partially supported by MURST.
Partially supported by DGICYT PB95-0354 and the Fundación del Amo, UCM.

must vanish on it with multiplicity odd on one hand and even on the other, so that the separation is actually impossible.

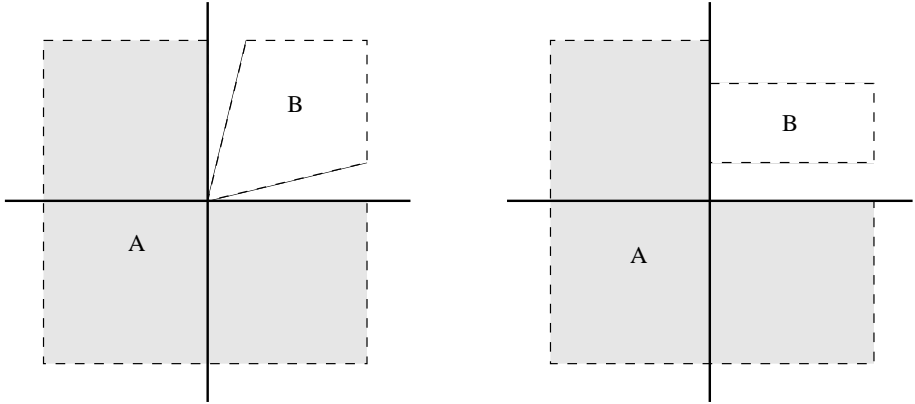

These easy examples contain, however, the main idea that inspired this paper, namely, that the relevant information for separation is at the boundary of $A$ and $B$ (including the points at infinity), although, as in the example above, this information may be hidden and appears only after a blowing-up, in a model where $A$ and $B$ are "fully" displayed.

In particular, we wonder whether there exists an "universal" obstruction to separation as the one shown in the examples, that is, the existence of a component of the boundary of $A$ and $B$ which is simultaneously *odd* and *even*. We show in Theorem 4.9 that this is indeed the case, provided that we look for it in the appropiate model of the variety. But allowing different models means that we are studying the problem "generically", i.e. up to codimension 1. Thus we divide our analysis of the separation procedure into two steps:

I: Study of the Generic Separation, that is, separation up to codimension one.
II: Analysis of Generic Separation versus Separation.

Generic separation is, obviously, a birational problem, that is, it depends only on the field $K$ of rational functions on $M$. In fact we have that $A$ and $B$ are separable if and only if their associated constructible sets $\widetilde{A}$ and $\widetilde{B}$ in the Space of Orderings $\Sigma_K$ of $K$ can be separated. Now, Bröcker's theorem cf. [Br3], [AnBrRz, Theorem IV.7.12], states that $\widetilde{A}$ and $\widetilde{B}$ cannot be separated if and only if there exists a finite subspace $X$ of $\Sigma_K$, in which $\widetilde{A}$ and $\widetilde{B}$ cannot be separated. Moreover it gives an upper bound for the chain length of such a subspace, so that we actually have an upper bound for the number of elements of $X$ and therefore for the number of functions on it. This way the problem is reduced, at least theoretically, to an infinite number of combinatorial ones, but still far away from being decidable.

We proceed as follows: our first step is to isolate a "nice" class of spaces of orderings which suffices to test the separation of $A$ and $B$. These are the *Geometric Spaces of Orderings* (GSO for short), which are introduced in Section 2. Roughly speaking they are spaces of orderings attached to discrete valuations and therefore they behave specially well with respect to the geometry of the variety, so that they allow to translate easily the "abstract" results into geometric statements. GSO are the natural generalization of the notion of algebroid fans considered in [AnRz2]. In Section 3 we prove the key result of the paper: the family of GSO is dense in the one of all spaces of orderings (Theorem 3.1). In particular this implies that GSO are enough to check the separation of $\widetilde{A}$ and $\widetilde{B}$ (Theorem 3.4). In Section 4 we use GSO to show that $A$ and $B$ are generically separated if and only if for any model of $M$, they are separated in an infinitesimal neighbourhood of their *walls* (the codimension 1 components of their boundaries), making precise the idea that obstacles to separation appear in the boundary, and showing the universal obstruction result mentioned above.

In Section 5 we show first that generic separation can be tested in a fixed model of $M$, provided that it is non–singular and the walls of $A$ and $B$ are at normal crossings. Then, under this assumptions we develop a geometric criterion (Theorem 5.5) which reduces the question of the generic separation of $A$ and $B$ in a neighbourhood of a wall $W$ to the separation of their *shadows* and *counter–shadows* in $W$. This way we lower the dimension by one and by induction we get a decision procedure for the generic separation (Theorem 6.1). This geometric criterion was proved for dimension 3 in [AcBgFo] by a different method.

Finally in Section 7 we detect the obstructions for two generically separable semialgebraic sets $A$ and $B$ to be separable, and show that they can be recursively tested, yielding the announced result on the decidability of the separation of semialgebraic sets of an affine variety over the real numbers (Theorem 7.6).

We want to thank A. Prestel who explained us the notion of recursive enumerability and the Model Theory aspects of the decision procedure during his visit to the University of Pisa. Also, this paper was completed while the second author enjoyed a stay at Stanford University partially supported by the "Fundación del Amo" of the Universidad Complutense of Madrid. He sincerely wants to thank both institutions for their help and warm atmosphere.

## 1. Notations and Preliminaries

The abstract theory of spaces of orderings was developed by Marshall in the series of papers [Mr1-4]. For a self-contained new presentation see [AnBrRz]. Here we just fix notations and definitions.

Let $K$ be a field, and consider its space of orderings $\Sigma_K = \mathrm{Spec}_r(K)$. Any element $\sigma \in \Sigma_K$, can be seen either as the cone in $K$ of all functions $f$ positive in $\sigma$, or as a signature $\sigma : K^*/\Sigma K^2 \to \{-1, 1\}$ which maps the class of the element $f \in K$ to 1 or $-1$ according to whether $f$ is positive or negative in $\sigma$. Thus, $\Sigma_K$ is a subset of the character group of $G_K = K^*/\Sigma K^2$, and the pair $(\Sigma_K, G_K)$, is an abstract Space of Orderings in the sense of Marshall. To keep a geometrical meaning in the notation we will write $f(\sigma) > 0$ instead of $\sigma(f) = 1$ and $f(\sigma) < 0$ instead of $\sigma(f) = -1$.

A *constructible* subset of $\Sigma_K$ is a set of the form
$$C = \bigcup_{i=1}^{p} \{\sigma \in \Sigma_K \mid f_{i1}(\sigma) > 0, ..., f_{ir_i}(\sigma) > 0\},$$
where $f_{ij} \in K$.

Given $E \subset K^*$ and $Y \subset \Sigma_K$, one defines
$$\begin{aligned} E^\perp &= \{\sigma \in \Sigma_K \mid f(\sigma) > 0, \text{ for all } f \in E\} \\ Y^\perp &= \{f \in K^* \mid f(\sigma) > 0, \text{ for all } \sigma \in Y\} = \bigcap_{\sigma \in Y} \sigma \end{aligned}$$

A subset $Y \subset \Sigma_K$ is called a *subspace* if it verifies $Y^{\perp\perp} = Y$. In this case the pair $(Y, G_K/Y^\perp)$ is also a subspace of orderings in the sense of Marshall.

Any singleton $E = \{\sigma\}$, $\sigma$ being an ordering, is a subspace which is called *atomic*. A set of 4 orderings $\{\sigma_1, \sigma_2, \sigma_3, \sigma_4\}$ with $\sigma_4 = \sigma_1\sigma_2\sigma_3$ is also a subspace which is called a 4–*element fan*. More generally a finite *fan* $F$ is a finite subset such that $F^3 \subset F$.

Two operations are defined between abstract spaces of orderings to build new ones: the *sum* $X_1 + X_2$ and the *extension* $X[H]$ by a 2-group $H$. Without entering in their definitions, let us say that in the context of the space $\Sigma_K$ of orderings of a field, the sum corresponds to the union of disjoint "independent" families $X_1$, $X_2$ of orderings of $K$, where independent means that $(X_1)^\perp \cdot (X_2)^\perp = K^*$, while the extension corresponds to the consideration of a family of orderings compatible with a real valuation and specializing to a given subspace of the residue field.

The structural theorem of finite spaces of orderings, [Mr1], [AnBrRz, Theorem IV.5.1], asserts that any finite subspace can be built in a unique way (up to isomorphism) by a finite sequence of sums and extensions, starting from a finite number of atomic spaces. Thus, any finite space of orderings, is represented, up to isomorphism, by a structure tree where the bottom points represent atomic spaces. Structure trees of finite spaces of orderings allow to work by induction along the tree, either to define properties or to prove them. For instance, that is the way we define the notion of Geometric Spaces of Orderings in the next section. Also the notions of stability index and chain length can be given in terms of the tree.

A valuation ring $V$ of $K$ is said *compatible* with an ordering $\sigma \in \Sigma_K$ if $\sigma$ makes convex the maximal ideal $\mathfrak{m}_V$ of $V$. Then $\sigma$ specializes to an ordering $\tau$ in the residue field $k_V$ of $V$, and we will write $\sigma \to \tau$. Conversely, given an ordering $\tau$ of $k_V$, Baer-Krull theorem ([BCR, Theorem 10.1.10]) asserts that there is a bijection between the set $F_\tau$ of orderings of $K$ compatible with $V$ and specializing to $\tau$, and the set of group homomorphisms $\phi : \Gamma/2\Gamma \to \{+1, -1\}$. More precisely, if $\sigma_0$ is an element of $F_\tau$, then

$$F_\tau = \{\sigma = (\phi, \sigma_0) \mid \phi \in \mathrm{Hom}(\Gamma/2\Gamma, \{+1, -1\})\}$$

where the sign of such a $\sigma$ on an element $f \in K$ is defined as $(\phi, \sigma_0)(f) = \phi(v(f))\sigma_0(f)$, $v$ being the valuation associated to $V$,

$$\sigma_0 F_\tau = \{(\phi, 1) \mid \phi \in \mathrm{Hom}(\Gamma/2\Gamma, \{+1, -1\})\}$$

can be identified with the group of characters $\mathrm{Hom}(\Gamma/2\Gamma, \{+1, -1\})$.

In particular, any ordering $\tau$ of the residue field can be "lifted" into $2^d$ different ways to $K$, where $d = \dim_{\mathbb{F}_2}(\Gamma/2\Gamma)$. A situation in which we will apply Baer-Krull theorem later, is the following:

**Construction 1.1** Let $B$ be a local regular ring with residue field $L$ and quotient field $K$. Suppose $\dim(B) = m$ and consider a system of parameters $x_1, \ldots, x_m$. By induction on $m$, we define a valuation $v_m$ in the quotient field $K$ of $B$ which dominates $B$, has residue field $L$ and value group $\mathbb{Z}^m$ ordered lexicographically.

If $m = 1$, then $B$ is a discrete valuation ring and we have the corresponding discrete rank one valuation $v_1$. For $m > 1$, we consider the discrete valuation ring $W = B_{(x_m)}$, whose valuation we denote by $w$. The residue field $K'$ of $W$ is the quotient field of the local regular ring $B' = B/(x_m)$. By induction we have in $K'$ a valuation $v_{m-1}$ with residue field $L$ and value group $\mathbb{Z}^{m-1}$. Then $v_m$ is the composite of $w$ and $v_{m-1}$. We denote its valuation ring by $V_m$. By construction $v_m(x_i) = (0, \ldots, 1, \ldots, 0)$ where the 1 appears in place $i$, and therefore $v_m(x_1), \ldots, v_m(x_m)$ are a basis of $\Gamma/2\Gamma = \mathbb{Z}^m/2\mathbb{Z}^m$, so that any homomorphism $\phi \in \operatorname{Hom}(\Gamma/2\Gamma, \{+1, -1\})$ may be identified with the mapping of $\{x_1, \ldots, x_m\}$ into $\{+1, -1\}$ which takes $x_i$ to $\phi(v_m(x_i))$.

Thus, if we fix an ordering $\tau$ in $L$ and denote by $F_\tau$ the set of orderings of $K$ compatible with $V_m$ and specializing to $\tau$, any element $\sigma \in F_\tau$ may be identified with the mapping of $\{x_1, \ldots, x_m\}$ into $\{+1, -1\}$ which takes $x_i$ to $\sigma(x_i)$.

Notice that if $\sigma_0$ is the (unique) generization of $\tau$ which is positive at all $x_i's$, then the map $\varphi : \sigma_0 F_\tau \longrightarrow \{1, -1\}^m$ defined by $\sigma_0 \sigma \longmapsto (\sigma(x_1), \ldots, \sigma(x_d))$, gives an isomorphism of $\mathbb{F}_2$-vector spaces. $\square$

We say that the valuation ring $V$ is *compatible with* a subspace of orderings $X \subset \Sigma$ if $V$ is compatible with every ordering $\sigma \in X$. It is easily checked that then the specializations of the orderings of $X$ form a subspace $X_V$ of $\operatorname{Spec}_r(k_V)$ called the *push-down* of $X$. Conversely, given a subspace of orderings $Y$ in $k_V$, the set $^V Y$ of orderings of $K$ wich are compatible with $V$ and specialize to an ordering of $Y$, is a subspace, called the *pull-back* of $Y$.

By Baer-Krull we have $^V Y = \widehat{H} \times Y$, where $\widehat{H} = \operatorname{Hom}(\Gamma/2\Gamma, \mathbb{F}_2)$, so that that $^V Y$ is the extension $^V Y = Y[\Gamma/2\Gamma]$. More generally, if $X \subset \Sigma_K$ is an extension, say $X = X'[H]$ where $X'$ is not an extension, let $V$ be the smallest valuation ring of $K$ compatible with $X$, and let $X_V$ be the push-down of $X$ and $^V X_V$ the pull-back of $X_V$. We have $X' = X_V$, and $X \subset {}^V X_V = X_V[\Gamma/2\Gamma]$, the inclusion being, in general, proper. In fact $H$ is isomorphic to the group $\Gamma/v(X^\perp)$, which in turn can be identified with a subgroup of $\Gamma/2\Gamma$, cf. [AnBrRz, Proposition IV.2.13]. Given an extension $X = X'[H]$, we will use the same notation introduced for the valuation rings, so that if $\sigma = (\widehat{h}, \sigma') \in X$, with $\sigma' \in X'$, $\widehat{h} \in \widehat{H}$, we will say that $\sigma$ *specializes* to $\sigma'$, and we will refer to the fiber $\widehat{H} \times \{\sigma'\}$ as the set of *generizations* of $\sigma'$.

## 2. Geometric Spaces of Orderings

From now on we assume that $K$ is a finitely generated extension of $\mathbb{R}$. In particular $K$ is the field of rational functions of some irreducible real algebraic variety $M$ and any such $M$ is called a *model* of $K$. It is well known that the topological dimension of any model $M$ coincides with the transcendence degree of $K$ over $\mathbb{R}$. Another useful fact is that we can always find *compact* models of $K$. This is immediate by taking the projective closure of any given model or considering its one-point compactification, which is always possible in the real case, [BCR], section 5.3.

**Definition 2.1.** *Let $K$ be a finitely generated extension of $\mathbb{R}$ of transcendence degree $n$. Let $(X, G)$ be a finite subspace of orderings of $\Sigma_K = \operatorname{Spec}_r(K)$. We define the notion of $(X, G)$ to be* geometric *by induction along its tree:*

1) *The atomic space $E = \{\sigma\}$ is geometric if its convex hull $W_\sigma$ is a discrete valuation of rank $n$.*
2) *If $X = X_1 + X_2$ then $X$ is geometric if and only if $X_1$ and $X_2$ are geometric.*
3) *If $X = X_1[H]$ is an extension with $\#(H) = 2^\alpha$ then $X$ is geometric if there is a discrete valuation ring $V$ of $K$ of rank $\alpha$ and with residue field $k_V$ a finitely generated extension of $\mathbb{R}$ of trnscendence degree $n - \alpha$, such that $X_1$ is a geometric subspace of $\operatorname{Spec}_r(k_V)$ and $X$ is the pull-back of $X_1$ by $V$.*

**Remark 2.2** Let us denote by $\mathfrak{B}_X$ the family of all valuation rings of $K$ compatible with some element of $X$. It follows from the definition that if $X$ is geometric then all the valuations of $\mathfrak{B}_X$ are discrete of maximum rank and with finitely generated residue field. In fact this can be taken as an alternative definition of GSO □

**Example 2.3** Let $M = \mathbb{R}^2$ and set $c(t) = (t, \xi(t))$, with $\xi(t) \in \mathbb{R}\{t\}$ a non algebraic power series (that is, $\xi(t)$ is not algebraic over the field of rational functions $\mathbb{R}(t)$). Then, the order $\sigma$ defined on $M$ by

$$\sigma(f) > 0 \quad \text{if } f(t, \xi(t)) = t^m u(t) \text{ with } u(0) > 0,$$
$$\sigma(f) < 0 \quad \text{if } f(t, \xi(t)) = t^m u(t) \text{ with } u(0) < 0.$$

is not geometric since $W_\sigma = \mathbb{R}[[t]] \cap \mathbb{R}(x,y)$ is of rank 1. On the contrary, if $\xi(t)$ is algebraic and $P(t,y)$ is its irreducible polynomial, the order $\sigma$ defined by the curve $P(x,y) = 0$ as,

$$\sigma(f) > 0 \quad \text{if } f = P^m Q \text{ with } Q \notin (P) \text{ and } Q(c(t)) > 0 \text{ for } t > 0,$$
$$\sigma(f) < 0 \quad \text{if } f = P^m Q \text{ with } Q \notin (P) \text{ and } Q(c(t)) < 0 \text{ for } t > 0.$$

is geometric. In this case $W_\sigma$ is the composite of $\mathbb{R}[x,y]_{(P)}$ and $(\mathbb{R}[x,y]/P)_a$, where $a = c(0)$. □

**Definition 2.4.** *We say that a finite geometric subspace of orderings $X$ is fully realized in an open semialgebraic subset $S$ of a compact model $M$ of $K$ if any $\sigma \in X$ specializes to a point in $S$ and its chain of specializations in $\mathrm{Spec}_r(\mathfrak{P}(M))$ has length $n$, i.e. there are $\sigma^{(i)} \in \mathrm{Spec}_r(\mathfrak{P}(M))$ such that $\sigma \to \sigma^{(1)} \to \cdots \to \sigma^{(n)} = \{x\} \in S$.*

As a first result on the abundance of geometric spaces of orderings, we have:

**Proposition 2.5.** *a) Any finite space of orderings of stability index $s$ is isomorphic to a geometric subspace realized in any open semialgebraic subset $S$ of dimension $s$ of a real algebraic variety $M$.*

*b) Let $X$ be a geometric space of orderings of $K$. Then there is a compact model of $K$ in which $X$ is fully realized. Moreover, for any compact model $M$ of $K$ there is sequence of blowings–up, $M_r \to M_{r-1} \to \cdots \to M_0 = M$, such that $X$ is realized in $M_r$.*

*Proof.* Part *a)* is [Br1, Proposition 3.3]. For part *b)* it is enough to blow-up $M$ till we reach a model in which all valuations (a finite number) of $\mathfrak{B}_X$ have a center of dimension equal to the transcendence degree of their residue fields. □

## 3. Density of geometric spaces of orderings

Fix an integer $k \geq 0$. Any subspace $X$ of $\Sigma_K$ with $k$ elements can be seen as a $k$-tuple in the product $\Sigma^k = \Sigma \times \cdots \times \Sigma$ ($k$ times). This identification is not bijective, unless we identify the tuples in $\Sigma^k$ up to permutations, but we will not care about this technicality, because it is irrelevant for our purposes. Anyway, the set $\Sigma^k$ carries the product topology of the Harrison topology of each factor space and we can discuss approximation properties in it. We recall that two finite subspaces of orderings are isomorphic if and only if they have the same construction tree, cf. §1.

In this section we will show the following:

**Theorem 3.1.** *Let $K$ be a finitely generated extension of $\mathbb{R}$ of transcendence degree $n$. Let $X \subset \Sigma_K$ be a subspace of orderings with $k$ elements. Then $X$ can be arbitrarily approximated in the Harrison topology by a geometric subspace $X'$ isomorphic to it.*

Before entering in the proof of the Theorem we need a couple of technical results:

**Lemma 3.2.** *Let $X \subset \Sigma$ be a subspace of orderings and assume that $X = X_1 + X_2$. If $X'_1$ and $X'_2$ are close enough to $X_1$ and $X_2$, then $X'_1 \cup X'_2 = X'_1 + X'_2$ as spaces of orderings.*

*Proof.* By the 4-element fan criterion, [Mr3, Corollary 7.5 *ii)*], [AnBrRz Theorem 7.2 *c)*], a union $X_1 \cup X_2$ of two disjoint closed subsets is a sum if and only if any 4–element fan is contained either in $X_1$ or in $X_2$, which in turn can be translated in terms that any element of $X_1$ (resp. $X_2$) can be separated from any group of three elements of $X_2$ (resp. $X_1$). Now, this is an open condition in the Harrison topology, so that the same is verified by any $X'_1$ and $X'_2$ close enough to $X_1$ and $X_2$, and the lemma follows at once. □

For the next result, given a field extension $K \subset L$ and a subspace of orderings $Y \subset \Sigma_L$ we denote by $Y \cap K$ the restriction of $Y$ to $K$, that is, $Y \cap K = \{\sigma \cap K \mid \sigma \in Y\}$.

**Lemma 3.3.** *Let $K \subset L$ be a field extension and let $Y \subset \Sigma_L$ be a finite subspace of orderings of $L$. Then there is a finitely generated subextension $K \subset F \subset L$ such that $X = Y \cap F$ is isomorphic to $Y$.*

*Proof.* We work by induction along the tree of $Y$. If $Y$ is atomic the result is trivial. If $Y = Y_1 + Y_2$, this means that $Y = Y_1 \cup Y_2$ and there is no 4-element fan not contained either in $Y_1$ or $Y_2$. By induction there is a finitely generated extension of $K$ such that $X_i = Y_i \cap F$ is isomorphic to $Y_i$. Moreover, adding if necessary, more elements to $F$, we may assume that $F$ contains functions which separates each element of $Y_1$ from each group of three of $Y_2$ and conversely. Thus, $X = X_1 \cup X_2 = X_1 + X_2$ and $X$ is isomorphic to $Y$.

Finally, let us assume that $Y$ is an extension, $Y = Y'[H]$, where $Y'$ is not an extension. Let $V$ be the smallest valuation ring of $L$ compatible with $Y$. Then $Y' \approx Y_V$, the push-down of $Y$, $H$ is isomorphic to a subgroup of $\Gamma_V/2\Gamma_V$ and $Y$ is a subspace of the pull-back $^V Y_V = Y_V[\Gamma_V/2\Gamma_V]$. Now, let $W = V \cap K$. $W$ is a valuation ring of $K$ with $k_W \subset k_V$ and $\Gamma_W \subset \Gamma_V$. By the induction hypothesis there is a finitely generated extension $k_1$ of $k_W$ with $k_1 \subset k_V$ and such that $X_1 = Y' \cap k_1$ is isomorphic to $Y'$. Let $a_1, \ldots, a_r \in V$ be such that their classes in $k_V$ do generate $k_1$ over $k_W$. Finally let $b_1, \ldots, b_s \in V$ be such that their values $v(b_1), \ldots, v(b_s)$ generate $H$ (remember that it is a finite group) and consider the extension $F = K(a_1, \ldots, a_r, b_1, \ldots, b_s)$. Then $V' = V \cap F$ is a valuation ring such that $\Gamma_{V'}/2\Gamma_{V'}$ contains $H$ as a subgroup and $k_1$ is a subfield of $k_V$ and it follows at once that $Y \cap F$ is isomorphic to $Y$. □

*Proof of Theorem 3.1:* Take a compact model $M$ of $K$. Then every polynomial is bounded on $M$, and it follows that every real valuation ring of $K$ contains the ring $\mathfrak{P}(M)$ of polynomial functions of $M$. Set $X = (\sigma_1, \ldots, \sigma_k)$ and let $U = U_1 \times \cdots \times U_k$ be an open neighbourhood of $X$ in $\Sigma^k$, with $U_i = \{f_{i1} > 0, \ldots, f_{ir_i} > 0\}$, $f_{ij} \in \mathfrak{P}(M)$. After shrinking the $U_i$'s we may assume that they are pairwise disjoint, and in particular the $f_{ij}$'s *separate* the orderings of $X$. We work by induction on the tree of $X$.

If $X = \{\sigma\}$ (the atomic space) then $U = \{f_1 > 0, \ldots, f_r > 0\}$. By Artin–Lang Theorem the semialgebraic set $S = \{f_1 > 0 \ldots, f_r > 0\} \cap \text{Reg}(M)$ is not empty. Take any point $x \in S$. Thus $\mathfrak{P}(M)_x$ is a regular local ring of dimension $n$ and by the Construction 1.2 there is a discrete valuation ring $V$ of $K$ of rank $n$, dominating the local ring $\mathfrak{P}(M)_x$. Finally take any ordering $\tau$ compatible with $V$. Obviously $\tau$ is geometric, $\tau \in U$ and $X \approx \{\tau\}$.

If $X = X_1 + X_2$, then, after reordering if necessary, $X_1 = (\sigma_1, \ldots, \sigma_\ell) \subset \Sigma^\ell$, $X_2 = (\sigma_{\ell+1}, \ldots, \sigma_k) \subset \Sigma^{k-\ell}$ and we have neighbourhoods $U^{(1)} = U_1 \times \cdots \times U_\ell$ and $U^{(2)} = U_{\ell+1} \times \cdots \times U_n$ of $X_1$ and $X_2$ respectively. Then by induction we find $X_1' \in U^{(1)}$ and $X_2' \in U^{(2)}$ isomorphic to them and geometric on $K$. Moreover, restricting the neighbourhoods if neccessary, by Lemma 3.2 above we have that $X' = X_1' \cup X_2' = X_1' + X_2'$. Thus, $X' \in U$ is geometric on $K$ and isomorphic to $X$, as wanted.

Next, let $X \approx X_1[H]$ be an extension with $\#(H) = 2^a$ and $X_1$ not an extension. Remember that by convention $2E$ is considered as the sum $E + E$ and not as an extension, so that $\#(X_1) \geq 2$. We have $X = \widehat{H} \times X_1$. Let $\pi : X \to X_1$ be the canonical projection and set $X_1 = \{\gamma_0, \ldots, \gamma_m\}$, $\widehat{H} = \{\hat{h}_0, \ldots, \hat{h}_{2^a-1}\}$. We order the elements of $X$ in a matrix arrangement:

$$\begin{array}{cccc} \sigma_0^0 & \sigma_0^1 & \cdots & \sigma_0^m \\ \sigma_1^0 & \sigma_1^1 & \cdots & \sigma_1^m \\ \vdots & \vdots & & \vdots \\ \sigma_{2^a-1}^0 & \sigma_{2^a-1}^1 & \cdots & \sigma_{2^a-1}^m \end{array}$$

where $\sigma_j^i = (\hat{h}_j, \gamma_i)$. Let $W$ be the minimal valuation ring compatible with $X$ and let $\Gamma$ be its value group. Let $X_W$ be the push-down of $X$ by $W$. Then $X_1 \approx X_W$, $X = X_1[H]$ is a subspace of $^W X_W = X_W[\Gamma/2\Gamma]$ and it follows at once that $r = \dim_{\mathbb{F}_2} \Gamma/2\Gamma \geq a$. From now on we identify $X_1$ and $X_W$. We also get that if we denote by $F_\gamma$ the pull-back of the ordering $\gamma \in \Sigma_{k_W}$ by $W$ then, for any $i = 0, \ldots, m$, we have $\{\sigma_1^i, \ldots, \sigma_{2^a}^i\} \subset F_{\gamma_i}$ (although this inclusion may be strict).

Now we apply resolution of singularities I and II ([Hk]), so that after finitely many blowings-up we may assume that $M$ is non-singular and all the $f_{ij}$'s are normal crossings. Let $\mathfrak{p} = \mathfrak{m}_W \cap \mathfrak{R}(M)$ be the center of $W$ in $\mathfrak{R}(M)$ where $\mathfrak{m}_W$ is the maximal ideal of $W$. Then $A = \mathfrak{R}(M)_\mathfrak{p}$ is a regular local ring of dimension, say, $d$, and has a regular system of parameters $x_1, \ldots, x_d$ such that, for all $i, j$

$$f_{ij} = u_{ij} x_1^{\alpha_{ij1}} \cdots x_d^{\alpha_{ijd}}$$

where the $u_{ij}$ are units of $A$ and the $\alpha_{ijk}$ are non-negative integers.

In this situation the quotient field $\kappa(\mathfrak{p})$ of $A/\mathfrak{p}$ is a subfield of the residue field $k_W$ of $W$ and performing some extra blowings-up we may assume also, by Lemma 3.3, that $\mathfrak{R}(M)$ contains enough functions so that the restrictions $\gamma_i \cap \kappa(\mathfrak{p})$ form, in $\text{Spec}_r(\kappa(\mathfrak{p}))$, a subspace isomorphic to $X_1$, which we still denote by $X_1$.

Next, notice that the signs of the elements $f_{ij}$ in any ordering $\sigma \to \gamma_i$ are completely determined by the signs of the parameters $x_l$ in $\sigma$ and the signs of the units (or more properly of their residue classes) in $\gamma_i$. Fix one of the $\gamma_i$, say $i = 0$. Thus, if we associate to each ordering $\sigma \in X$ the $d$-tuple $(\sigma(x_1), \ldots, \sigma(x_d))$, these $d$-tuples must differenciate the orderings over $\gamma_0$. In particular we get that $2^d \geq 2^a$, i.e. $a \leq d$. We will see that after some operations which in the end will be quadratic transforms of $A$ (and hence of $M$) we will find a system of parameters $x_1, \ldots, x_a, x_{a+1}, \ldots, x_d$ such that $X$ will coincide with the product $\widehat{G} \times X_1$, where $\widehat{G}$ is the subspace of the mappings $\varepsilon$ of $\{x_1, \ldots, x_d\}$ into $\{-1, +1\}$ such that $\varepsilon(x_\ell) = +1$ for all $\ell > a$.

To make this precise let us come back to our matricial distribution of the elements of $X$. Take any ordering over $\gamma_0$ say $\sigma_0^0$. Changing $x_j$ by $-x_j$ if necessary, we may assume that $\sigma_0^0(x_j) = 1$ for all $j$ (this $\sigma_0^0$ will define the zero homomorphism of $\{x_1, \ldots, x_d\}$ into $\mathbb{F}_2 = \{-1, +1\}$). Now take $a$ more elements over $\gamma_0$, say $\sigma_1^0, \ldots, \sigma_a^0$ which are linearly independent (and therefore $\sigma_0^0, \sigma_1^0, \ldots \sigma_a^0$ generate the fiber over $\gamma_0$). Moreover, take now an element in each of the other columns, say $\sigma_0^1, \ldots, \sigma_0^m$. Notice that the elements $\sigma_0^0, \ldots, \sigma_a^0, \sigma_0^1, \ldots, \sigma_0^m$, generate $X$, since they generate the first column and with this and the selected elements in each of the other columns we can get any element of $X$ by completing squares of 4-element fans: $\sigma_j^i = \sigma_j^0 \sigma_0^i \sigma_0^0$.

**Claim:** *After some additional blowing-ups, we find a regular local ring $B$ dominating $A$, with the same residue field and a system of parameters $y_1, \ldots, y_d$ of $B$ such that all $f_{ij}$'s are normal crossings in $B$ with respect to them and for all $j = 0, \ldots, a$ it holds*

$$\sigma_j^0(y_k) = \begin{cases} +1 & \text{for} \quad k \geq j+1 \\ -1 & \text{if} \quad k = j \end{cases}$$

In other words, we are saying that the matrix of signs of the $\sigma_j^0$, $j = 1, \ldots, a$, at the system of parameters is low triangular with $-1$ in the main diagonal, i.e. we have an arrangement of the type:

|  | $y_1$ | $y_2$ | $y_3$ | $\cdots$ | $y_{a-1}$ | $y_a$ | $y_{a+1}$ | $y_{a+2}$ | $\cdots$ | $y_d$ |
|---|---|---|---|---|---|---|---|---|---|---|
| $\sigma_0^0$ | 1 | 1 | 1 | $\cdots$ | 1 | 1 | 1 | 1 | $\cdots$ | 1 |
| $\sigma_1^0$ | $-1$ | 1 | 1 | $\cdots$ | 1 | 1 | 1 | 1 | $\cdots$ | 1 |
| $\sigma_2^0$ | $*$ | $-1$ | 1 | $\cdots$ | 1 | 1 | 1 | 1 | $\cdots$ | 1 |
| $\vdots$ | $\vdots$ | $\vdots$ | $\vdots$ |  | $\vdots$ | $\vdots$ | $\vdots$ | $\vdots$ |  | $\vdots$ |
| $\sigma_{a-1}^0$ | $*$ | $*$ | $*$ | $\cdots$ | $-1$ | 1 | 1 | 1 | $\cdots$ | 1 |
| $\sigma_a^0$ | $*$ | $*$ | $*$ | $\cdots$ | $*$ | $-1$ | 1 | 1 | $\cdots$ | 1 |

where $*$ represents $+1$ o $-1$.

In fact, taking into account that the functions $f_{ij}$ separate the orderings of $F_{\gamma_0}$ and that $\sigma_0^0, \sigma_1^0, \ldots, \sigma_a^0$ generate this set, they must go to independent vectors under the map

$$\varphi : \sigma_0^0 F_{\gamma_0} \longrightarrow \{1, -1\}^m$$

defined by

$$\sigma_0^0 \sigma \longmapsto (\sigma(x_1), \ldots, \sigma(x_d)),$$

Now, the proof of the claim consists of translating to our situation the classical Gauss elimination algorithm to make low-triangular the matrix of signs of the images of the $\sigma$'s by these mappings, using the pivot method of beginners Linear Algebra course. Remember that in our case the zero element of $\mathbb{Z}_2$ is 1 while the non-zero element is $-1$ and that taking the difference between columns $k$ and $\ell$ of the matrix corresponds, with our multiplicative notation for $\mathbb{Z}_2$ to take the quotient of the column $k$ over column $\ell$. This, in turn, yields to replace the parameter $x_k$ by $x_k/x_\ell$, which corresponds to perform a blowing-up in $A$.

Let us see it formally. Since $\sigma_1^0$ cannot have the same signs as $\sigma_0^0$ on all the parameters (i.e. $\varphi(\sigma_1^0)$ is not the zero vector), there is some $k$ such that $\sigma_1^0(x_k) \neq \sigma_0^0(x_k)$. We reorder the parameters so that $\sigma_0^1(x_l) = -1$ for $l < r$ and $\sigma_0^1(x_l) = 1$ for $l \geq r$. Consider the extension
$$A^{(1)} = A[x_2/x_1, \ldots, x_{r-1}/x_1]_{(x_1, x_2/x_1, \ldots, x_{r-1}/x_1, x_r, \ldots, x_d)}.$$
We set $x_1^{(1)} = x_1$, $x_k^{(1)} = x_k/x_1$ for $2 \leq k \leq r-1$ and $x_k^{(1)} = x_k$ for $k \geq r$. Now $A^{(1)}$ is a local regular ring dominating $A$, with the same residue field and $x_1^{(1)}, x_2^{(1)}, \ldots, x_d^{(1)}$ is a regular system of parameters of $A^{(1)}$. Furthermore the expression
$$f_{ij} = u_{ij} x_1^{\alpha_{ij1}} \cdots x_d^{\alpha_{ijd}}$$
can also be written
$$f_{ij} = u_{ij} (x_1^{(1)})^{\alpha_{ij1} + \alpha_{ij2} + \cdots + \alpha_{ij(r-1)}} (x_2^{(1)})^{\alpha_{ij2}} \cdots (x_d^{(1)})^{\alpha_{ijd}}.$$
which means that the $f_{ij}$ are still normal crossings in $A^{(1)}$, so that all conditions verified by $A$ are similarly verified by $A^{(1)}$. Moreover, we have $\sigma_1^0(x_1^{(1)}) = -1$ and $\sigma_1^0(x_k^{(1)}) = 1$ for all $k \geq 2$, so that we have constructed the first step in the induction proccess.

Assume that we have already found a local regular ring $A^{(\ell)}$ dominating $A$ with the same residue field that the latter and with a system of parameters $x_1^{(\ell)}, \ldots, x_d^{(\ell)}$ such that the $f_{ij}$ are normal crossings for them in $A^{(\ell)}$ and for all $0 \leq j \leq \ell$ it holds $\sigma_j^0(x_k^{(\ell)}) = 1$ for all $k \geq j+1$ and $\sigma_j^0(x_j^{(\ell)}) = -1$. We construct $A^{(\ell+1)}$ as follows:

Consider $\sigma_{\ell+1}^0$. We claim that there is $k \geq \ell+1$ such that $\sigma_{\ell+1}^0(x_k^{(\ell)}) = -1$. For otherwise, a look at the above table of signs shows at once that $\sigma_0^0 \sigma_{\ell+1}^0$ would be in the subspace generated by $\sigma_0^0 \sigma_1^0, \ldots, \sigma_0^0 \sigma_\ell^0$, against our assumption that $\sigma_0^0, \ldots \sigma_a^0$ are affine independent. Then, after reordering $x_{\ell+1}^{(\ell)}, \ldots, x_d^{(\ell)}$, we may assume that $\sigma_{\ell+1}^0(x_k^{(\ell)}) = 1$ for $k \geq r$ and $\sigma_{\ell+1}^0(x_k^{(\ell)}) = -1$ for $\ell+1 \leq k < r$. Consider the extension
$$A^{(\ell+1)} = A^{(\ell)}[x_{\ell+2}^{(\ell)}/x_{\ell+1}^{(\ell)}, \ldots, x_{r-1}^{(\ell)}/x_{\ell+1}^{(\ell)}]_{(x_1^{(\ell)}, \ldots, x_{\ell+1}^{(\ell)}, x_{\ell+2}^{(\ell)}/x_{\ell+1}^{(\ell)}, \ldots, x_{r-1}^{(\ell)}/x_{\ell+1}^{(\ell)}, x_r^{(\ell)}, \ldots, x_d^{(\ell)})},$$
and set $x_k^{(\ell+1)} = x_k^{(\ell)}$ for $1 \leq k \leq \ell+1$, $x_k^{(\ell+1)} = x_k^{(\ell)}/x_{\ell+1}^{(\ell)}$ for $\ell+2 \leq k \leq r-1$ and $x_k^{(\ell+1)} = x_k^{(\ell)}$ for $k \geq r$. An immediate computation shows that for $1 \leq j \leq \ell+1$ it holds $\sigma_j^0(x_k^{(\ell+1)}) = 1$ for all $k \geq j+1$ and $\sigma_j^0(x_j^{(\ell+1)}) = -1$, so that we have done the step $\ell+1$. By induction, this shows the claim.

Now, let us go back to the matrix representation of $X$ and notice that the signs of the parameters $y_{a+1}, \ldots, y_d$ produce a partition of $X$ which is compatible with the columns. In fact, suppose that $\sigma_0^i(y_k) = \varepsilon$ for some $k = a+1, \ldots, d$, where $\varepsilon = 1, -1$. Then for any other $j = 1, \ldots, 2^a - 1$ we have $\sigma_j^i(y_k) = \sigma_0^0(y_k)\sigma_0^i(y_k)\sigma_j^0(y_k) = \sigma_0^i(y_k) = \varepsilon$, in other words all members of the column $i$ have the same sign at any of the $y_k$ for $k = a+1, \ldots, d$. In particular, the family $S_\varepsilon = \{\varepsilon_{a+1} y_{a+1} > 0, \ldots, \varepsilon_d y_d > 0\}$, where $\varepsilon = (\varepsilon_{a+1}, \ldots, \varepsilon_d) \in \{-1, 1\}^{d-a}$ produces a partition of $X$ compatible with the columns.

Taking into account the structure of $X = X_1[H]$ as extension, this means that for each $k > a$ the function $y_k$ defines a function over $X_1$. Setting $X_{1,\varepsilon} = X_1 \cap S_\varepsilon$ we have a partition $X_1 = \bigcup_\varepsilon \{X_{1,\varepsilon}\}$ of $X_1$ into disjoint basic sets.

Now we consider the diagram

$$\begin{array}{ccc} B & \longrightarrow & C = B_{(y_1,\ldots,y_a)} \subset K \\ \downarrow & & \downarrow \\ B/(y_1,\ldots,y_a) & \longrightarrow & k_C \\ \downarrow & & \\ k_A = k_B & & \end{array}$$

where $k_B$, $k_C$ stand for the residue fields of $B$ and $C$ respectively. By construction, these two residue fields are finitely generated over $\mathbb{R}$. Now let $k'_B$ be a quasicoefficient field of $B$, that is a subfield $k'_B \subset B$, such that the extension $k'_B \subset k_B$ induced by the canonical homomorphism $B \to k_B$ is algebraic (even finite in our case). Then, since $k_C$ is the quotient field of the ring $B/(y_1,\ldots,y_a)$, which is local regular of dimension $d - a$ and $k_B = k_A = \kappa(\mathfrak{p})$, we get

$$\mathrm{tr.deg.}[k_C : \mathbb{R}] = \mathrm{tr.deg.}[k_C : k'_B] + \mathrm{tr.deg.}[k'_B : \mathbb{R}] \geq$$
$$d - a + \mathrm{tr.deg.}[\kappa(\mathfrak{p}) : \mathbb{R}] = d - a + \dim(\mathfrak{R}(M)/\mathfrak{p}) =$$
$$d - a + \dim(\mathfrak{R}(M)) - \mathrm{ht}(\mathfrak{p}) = d - a + n - d = n - a.$$

Let us lift our space $X_1$ from $k_B = \kappa(\mathfrak{p})$ to $k_C$. Take $\gamma_i \in X_1$, and assume that $\gamma_i \in S_\varepsilon$, i.e., $\gamma_i(y_k) = \varepsilon_k$ for any $k = a+1,\ldots,d$. Since $B/(y_1,\ldots,y_a)$ is a local regular ring with parameters $y_{a+1},\ldots,y_d$, we can lift $\gamma_i$ to a unique ordering $\tau_i$ of $k_C$ such that $\gamma_i(y_j) = \tau_i(y_j) = \varepsilon_k$ for $a+1 \leq j \leq d$.

We claim that $Y_1 = \{\tau_1,\ldots,\tau_m\}$ is a subspace of orderings of the extension $\widehat{G} \times X_1$ where $\widehat{G}$ is the group of mappings of $\{y_{a+1},\ldots,y_d\}$ into $\{+1,-1\}$, and that $Y_1 \approx X_1$. Indeed, notice that as element of $\widehat{G} \times X_1$ we have $\tau_i = (h, \gamma_i)$, where $h$ is the function such that $h(y_k) = \varepsilon_k = \gamma_i(y_k)$ for $k = a+1,\ldots,d$. To check that it is a subspace we have to see that for any four element fan $F \subset \widehat{G} \times X_1$, $\#(F \cap Y_1) \neq 3$, cf. [AnBrRz, Corollary IV.1.8]. Assume that $F = \{(g_i, \alpha_i)\}_{1 \leq i \leq 4}$, with, say $(g_1, \alpha_1), (g_2, \alpha_2), (g_3, \alpha_3) \in Y_1$. Then $F_1 = \{\alpha_i\}_{1 \leq i \leq 4}$ is a fan of $X_1$, that is, $\alpha_1 \alpha_2 \alpha_3 = \alpha_4$. Thus,

$$g_4(y_q) = g_1(y_q) g_2(y_q) g_3(y_q) = \alpha_1(y_q) \alpha_2(y_q) \alpha_3(y_q) = \alpha_4(y_q)$$

what means that $(g_4, \alpha_4) \in Y_1$ as wanted. That $Y_1 \approx X_1$ is now immediate since the projection from the former to the second is onto with fibers consisting of a single point.

To finish, by induction, let $X'_1$ be a geometric subspace of $\mathrm{Spec}_r(k_C)$ isomorphic to $Y_1$ and close enough to it so that each element $\gamma'_i \in X'_1$ verifies that $\gamma'_i(y_k) = \tau_i(y_k) = \gamma_i(y_k)$ for all $k = a+1,\ldots,d$ and also $\gamma'_i(u_{rt}) = \tau_i(u_{rt}) = \gamma_i(u_{rt})$ on the units $u_{rt}$ coming from the functions $f_{rt}$.

Since $C$ is local regular with parameters $y_1,\ldots,y_a$, we can produce an extension $X' = \widehat{H} \times X'_1$ with $\widehat{H}$ equal to the set of mappings of $\{y_1,\ldots,y_a\}$ into $\{+1,-1\}$ using a discrete valuation ring of rank $a$ which dominates $C$, as in the Construction 1.2. By construction $X'$ is geometric and $X' \approx X$. Moreover, if we set $\beta^i_j = (h_j, \gamma'_i)$, where $h_j \in \widehat{H}$ is the function defined by $h_j(y_k) = \sigma^0_j(y_k)$ for $1 \leq k \leq a$, we have that $\beta^i_j(f_{rt}) = \sigma^i_j(f_{rt})$ for all $i,j,r,t$, and since the $f_{rt}$'s define the neighborhood $U \subset \Sigma_K$ of $X$ fixed at the beginning, we conclude $X' \subset U$, which completes the proof.                                                                                                    □

As a consequence we have the following result. Remember that two constructible subsets $\widetilde{A}$ and $\widetilde{B}$ of $\Sigma_K$ are called separable if there exists $g \in K$ such that $\widetilde{A} \subset \{g > 0\}$ and $\widetilde{B} \subset \{g < 0\}$:

**Theorem 3.4.** *Let $K$ be a finitely generated extension of $\mathbb{R}$ of transcendence degree $n$. Let $\widetilde{A}$ and $\widetilde{B}$ be two disjoint constructible subsets of $\Sigma_K$. Then $\widetilde{A}$ and $\widetilde{B}$ are separable if and only if for any finite GSO $Y \subset \Sigma_M$ of chain-length $\leq 2^{n-1}$, $\widetilde{A} \cap Y$ and $\widetilde{B} \cap Y$ are separable in $Y$.*

*Proof.* One direction is trivial after Bröcker's theorem, [Br2],[AnBrRz Theorem IV.7.12]. For the other, suppose that $\widetilde{A}$ and $\widetilde{B}$ are not separable. Then, again by Bröcker's theorem there exists a finite subspace $Y$ of $\Sigma_M$ with chain-length $\leq 2^{n-1}$ such that $\widetilde{A} \cap Y$ and $\widetilde{B} \cap Y$ are not separable. Suppose that $\#(Y) = k$ and order its elements as a $k$-tuple $(\sigma_1, \ldots, \sigma_k)$. Also $Y$ has a finite number of different functions, namely the elements of the group $G_Y = K/Y^\perp$. Take $f_1, \ldots, f_N \in K$ such that
$$G_Y = \{1Y^\perp, -1Y^\perp, f_1 Y^\perp, \ldots, f_N Y^\perp\},$$
and let $f_{N+1}, \ldots, f_r$ be a family of elements of $K$ which describes the sets $A$ and $B$. Then, for any $j = 1, \ldots, k$ consider the neigbourhood of $\sigma_j$:
$$U_j = \{f_i(\sigma_j) f_i > 0 \mid i = 1, \ldots r\}$$
and set $U = U_1 \times \cdots \times U_k$. Thus $U$ is a neighbourhood of $Y$ in $\Sigma_K^k$ and by the density Theorem 3.1, there is a GSO, $X = (\tau_1, \ldots, \tau_k) \in U$, isomorphic to $Y$.

First of all notice that by our construction $\sigma_j \in A$ if and only if $\tau_j \in A$ and
$$G_X = \{1X^\perp, -1X^\perp, f_1 X^\perp, \ldots, f_N X^\perp\}.$$
Indeed, since $X$ and $Y$ are isomorphic we have $\#(G_X) = \#(G_Y)$. On the other hand, by construction all $f_j$ are non trivial and different on $X$, so that they fill the whole $G_X$. But now, since the $f_j$'s have the same signs over the $\sigma_s$'s and the $\tau_s$'s we have that $A \cap X$ and $B \cap X$ are not separable.                □

## 4. Generic Separation and Walls

Now let $M$ be an irreducible algebraic variety over $\mathbb{R}$ (which we assume compact), let $\mathfrak{R}(M)$ be the ring of regular functions on $M$ and $K(M)$ its quotient field. Let $A, B$ be disjoint semialgebraic subsets of $M$.

**Definition 4.1.**   a) *We say that $A$ and $B$ are* separable *if there exists a regular function $f \in \mathfrak{R}(M)$ such that $f(A) > 0$ and $f(B) < 0$.*
  b) *We say that $A$ and $B$ are* generically separable *if there exists a non-zero regular function $f \in \mathfrak{R}(M)$ such that $f(A) \geq 0$ and $f(B) \leq 0$. Equivalently, $A$ and $B$ are generically separable if there exist a proper algebraic subset $Y \subset M$ such that $A \setminus Y$ and $B \setminus Y$ are separable.*

Now let us denote by $\Sigma_{\mathfrak{R}(M)}$ or simply by $\Sigma_M$ the space of orderings of $\mathfrak{R}(M)$. Remember that we have a tilde map which assigns to any semialgebraic set $S \subset M$ the constructible subset $\widetilde{S} \subset \Sigma_M$ defined by the same equations as $S$. Two semialgebraic sets have the same tilde image if and only if they are generically equal, i.e. they are equal up to a subset of codimension at least one, [BCR, Proposition 7.6.3].

Thus the study of constructible subsets of $\Sigma_M$ translates into generic properties of semialgebraic sets, that is, what happens up to a set of smaller dimension. In particular we have that $A$ and $B$ are generically separable if and only if $\widetilde{A}$ and $\widetilde{B}$ are separable in $\Sigma_M$.

From now till Section 7 we will be dealing with the generic separation of $A$ and $B$. Since this is a birational property, we may assume that *M is non-singular and that A and B are open*. We are interested in testing the generic separation of $A$ and $B$ in terms of their boundaries. To start with, recall that if $\overline{A} \cap \overline{B} = \emptyset$ then $A$ and $B$ are separable (since $M$ is compact). So, only the case in which $A$ and $B$ have some common boundary is interesting.

**Definition 4.2.** *We define the* walls *of A and B in M as the irreducible $(n-1)$-dimensional components of the Zariski closure of the boundaries of A and B, that is, the irreducible components of $\overline{\partial A}^z \cup \overline{\partial B}^z$ of dimension $n-1$.*

Notice that since we are assuming that $M$ is non-singular, for any wall $W$, the localization $\mathfrak{R}(M)_{\mathfrak{J}(W)}$ of the ring of regular functions on $M$ at the prime ideal defined by $W$ is a rank 1, discrete, valuation ring of $K$ with residue field $\mathfrak{K}(W)$. Hence, the pull-back of the space of orderings, $\Sigma_W$, of $\mathfrak{K}(W)$ by the valuation $\mathfrak{R}(M)_{\mathfrak{J}(W)}$ is a subspace of orderings of $\Sigma_K$. Moreover, this pull-back is isomorphic to the extension $\Sigma_W[\mathbb{Z}_2]$. Thus, for any wall we see the extension $\Sigma_W[\mathbb{Z}_2]$ as a subspace of $\Sigma_K$. From now on we fix the notation $\widehat{\mathbb{Z}}_2 = \{1, i\}$, where 1 stands for the trivial map from $\mathbb{Z}_2$ to $\mathbb{Z}_2$ and $i$ for the identity. In particular for any $\gamma \in \Sigma_W$ we have two generizations $\alpha$ and $i\alpha$ in $\Sigma_W[\mathbb{Z}_2]$. Also, given a subset $\widetilde{C} \subset \Sigma_W[\mathbb{Z}_2]$ we denote by $\widetilde{C}_W$ the set of its specializations in $\Sigma_W$.

From a geometric point of view, the space $\Sigma_W[\mathbb{Z}_2]$ represents a kind of infinitesimal neighbourhood of $W$ in $M$, and the basic idea behind scenes is that the $A$ and $B$ can be separated if and only if they can be separated in such a neighbourhoods (compare with [AcBgFo, §1]). Indeed, as a first result we have:

**Proposition 4.3.** *A and B are generically separable if and only if for any non-singular compact model M of K and any wall of A and B in M, $\widetilde{A}$ and $\widetilde{B}$ are separable in $\Sigma_W[\mathbb{Z}_2]$.*

*Proof.* If $A$ and $B$ are generically separable, $\widetilde{A}$ and $\widetilde{B}$ are separable in $\Sigma_K$, and therefore so are they in any subspace $\Sigma_W[\mathbb{Z}_2]$, which shows the trivial direction.

For the other, assume that $A$ and $B$ are not generically separable. Then there exists a finite subspace of orderings $X \subset \Sigma_K$ such that $\widetilde{A} \cap X$ and $\widetilde{B} \cap X$ are not separable. Take $X$ with minimal cardinal. By Theorem 3.1 we may assume that $X$ is a geometric space of orderings. Furthermore, since for a sum $X = X_1 + X_2$ we have that $\widetilde{A} \cap X$ and $\widetilde{B} \cap X$ are separable if and only if $\widetilde{A} \cap X_i$ and $\widetilde{B} \cap X_i$ are separable for $i = 1, 2$, by the minimality assumption we have that $X$ is an extension, say $X = Y[\mathbb{Z}_2]$, with $Y$ a geometric space of orderings. Now, take a model $N$ of $K$ in which $X$ is fully realized, cf. Proposition 2.5 b). Then $Y$ is a GSO of a hypersurface $W$ of $N$ and we claim that $W$ is a wall of $A$ and $B$ in $N$.

Indeed, suppose that $W \not\subset \overline{\partial A}^z \cup \overline{\partial B}^z$. Since we are assuming that $\widetilde{A}$ and $\widetilde{B}$ are not separated in $X$, we have that $\widetilde{A} \cap X \neq \emptyset$, $\widetilde{B} \cap X \neq \emptyset$ and since any $\sigma \in X$ specializes to an ordering of $Y$, we get in particular that $W \cap \overline{A} \neq \emptyset$ and $W \cap \overline{B} \neq \emptyset$. Now, since $W$ is not in the boundary of $A$, for any $\alpha \in \widetilde{A} \cap X$ we also have $i\alpha \in \widetilde{A} \cap X$ and the same happens for $B$. Therefore $\widetilde{A} \cap X$ and $\widetilde{B} \cap X$ are saturated for the specialization map $X \to Y$. This shows that $\widetilde{A} \cap X$ and $\widetilde{B} \cap X$ are separable if and only if their specializations $\widetilde{A}_W$ and $\widetilde{B}_W$ are separable in $Y$, so that, in particular these are not separable. Finally consider the subspace $\{1\} \times Y \subset X$. The specialization map $X \to Y$ defines an isomorphism from $\{1\} \times Y$ onto $Y$ sending $\widetilde{A} \cap (\{1\} \times Y)$ to $\widetilde{A}_W \cap Y$ and $\widetilde{B} \cap (\{1\} \times Y)$ to $\widetilde{B}_W \cap Y$. In particular this shows that $\widetilde{A} \cap (\{1\} \times Y)$ and $\widetilde{B} \cap (\{1\} \times Y)$ are not separable, against our assumption on the minimality of $X$. □

**Remark 4.4** The consideration of all models is important. For instance, going back to the example considered in the second picture of the introduction we see that $A$ and $B$ are not generically separable, but they are separable in all spaces $\Sigma_W[\mathbb{Z}_2]$ for the walls $W$ in the picture of the left. The wall that produces the obstruction is the $y$–axis of the picture on the right, which only appears after blowing-up. □

Looking for a more geometrical description of the obstruction to the separation we introduce the following definition, where $M$, $A$ and $B$ are as above:

**Definition 4.5.** *We say that a wall $W$ is* odd *if there are $A' \subset A$, $B' \subset B$ which are generically separable and such that any function $f$ separating $A'$ and $B'$ changes sign across $W$. Similarly, we say that a wall $W$ is* even *if there are $A' \subset A$, $B' \subset B$ which are generically separable and such that any function $f$ separating $A'$ and $B'$ does not change sign across $W$.*

Intuitively a wall is odd if any hypothetic function $f$ which would separate $A$ and $B$ must change sign along $W$ and even if it can not change sign along $W$. In algebraic terms this means that if $t$ is a uniformizer of the local ring $\mathfrak{R}(M)_{\mathfrak{J}(W)}$, then, in this ring we have $f = t^m u$ where $u$ is a unit and $m$ is odd (resp. even) if and only if $W$ is odd (resp. even). Here are some examples of odd and even walls (cf. also [AcAnBg]):

**Example 4.6** The wall $W$ is odd in the pictures:

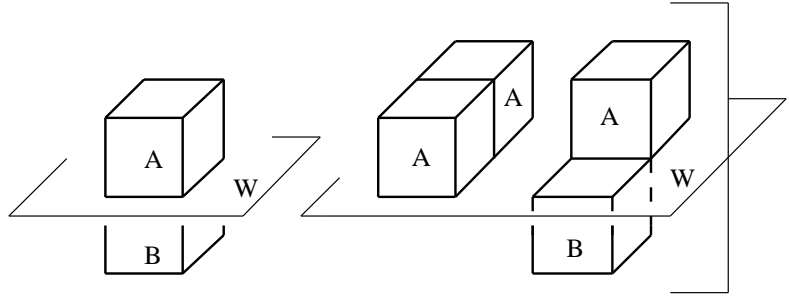

where $A' = A$ and $B' = B$, and is even in

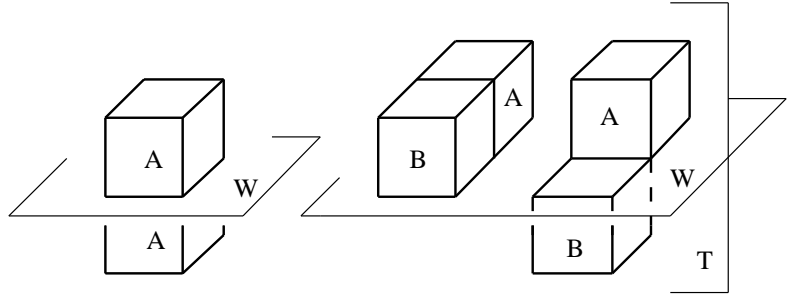

where we take also $A' = A$ and $B' = B$ and the plane $T$ as separating function in the second case.

Obviously, if some wall $W$ is simultaneously odd and even then $A$ and $B$ are not generically separable. Here are two standard examples of this situation,

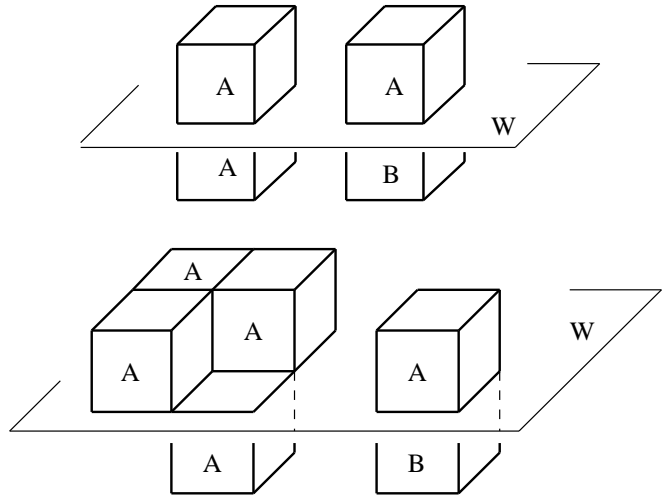

where in the second, the four cubes $A$ on the left cause $W$ to be even, while the couple $A$, $B$ on the right make $W$ odd. □

The amazing thing is that the converse is also true: $A$ and $B$ are not generically separable if and only if there is a wall which is simultaneously odd and even, provided that we look at the walls in the appropiated model. Before explaining this, we have to make a new small detour into spaces of orderings.

The notions of odd and even walls have a nice translation in terms of the space of orderings $\Sigma_W[\mathbb{Z}_2]$ which in fact helps often to recognize whether $W$ is odd or even. We define the *saturations* $\widetilde{A}^{\#}$ and $\widetilde{B}^{\#}$ of $\widetilde{A}$ and $\widetilde{B}$ in $\Sigma_M$ as:

$$\widetilde{A}^{\#} = \{\sigma \in \Sigma_K \mid \sigma \text{ can be written as a product of}$$
$$\text{elements of } \widetilde{A} \cup \widetilde{B} \text{ with an odd number of them in } \widetilde{A}\}$$

$$\widetilde{B}^{\#} = \{\sigma \in \Sigma_K \mid \sigma \text{ can be written as a product of}$$
$$\text{elements of } \widetilde{A} \cup \widetilde{B} \text{ with an odd number of them in } \widetilde{B}\}$$

We remember that since $\sigma$ is an ordering, it must always be a product of an odd number of elements of $A \cup B$. Intuitively we are considering all the orderings such that if $f$ is a function with $f|_{\widetilde{A}} > 0$ and $f|_{\widetilde{B}} < 0$, then also $f|_{\widetilde{A}^\#} > 0$ and $f|_{\widetilde{B}^\#} < 0$. This will be completely clear after the following result which collects the basic idea behind Bröcker's separation theorem [Br2], [AnBrRz, Theorem IV.7.12]:

**Proposition 4.7.** *The following are equivalent:*
  a) $\widetilde{A}$ *and* $\widetilde{B}$ *are separable;*
  b) *Any* $A'$, $B'$ *with* $\widetilde{A} \subset A' \subset \widetilde{A}^\#$ *and* $\widetilde{B} \subset B' \subset \widetilde{B}^\#$ *are separable.*
  c) $\widetilde{A}^\# \cap \widetilde{B}^\# = \emptyset$.

*Proof.* Lets us see that *a)* implies *c)*. Suppose that $f$ is a function with $f|_{\widetilde{A}} > 0$ and $f|_{\widetilde{B}} < 0$ and assume that there is $\sigma \in \widetilde{A}^\# \cap \widetilde{B}^\#$. Then, from $\sigma \in \widetilde{A}^\#$ we get that $\sigma$ is a product with an odd number of terms in $\widetilde{A}$ and an even number in $\widetilde{B}$. Hence $f(\sigma) > 0$. On the other hand, similarly, from $\sigma \in \widetilde{B}^\#$ we get that $f(\sigma) < 0$, contradiction.

Next, to see that *c)* implies *b)* it is enough to show that if $\widetilde{A}^\# \cap \widetilde{B}^\# = \emptyset$ then $\widetilde{A}^\#$ and $\widetilde{B}^\#$ are separable. Let $X$ be the subspace of $\Sigma_K$ spanned by $\widetilde{A} \cup \widetilde{B}$. Then we have $X = (\widetilde{A})^\# \cup (\widetilde{B})^\#$. Moreover, by [AnBrRz, Theorem IV.7.11] there is a constant $\rho$ such that any element spanned by $\widetilde{A} \cup \widetilde{B}$ is a product of at most $\rho$ elements of $\widetilde{A} \cup \widetilde{B}$. This implies in particular that $\widetilde{A}^\#$ and $\widetilde{B}^\#$ are closed, and since they are disjoint also open in $X$. Now, take any 4-element fan $\{\alpha_1, \alpha_2, \alpha_3, \alpha_4\}$ in $X$ with $\alpha_1, \alpha_2, \alpha_3 \in \widetilde{A}^\#$. Then by the very definition of $\widetilde{A}^\#$ we get $\alpha_4 = \alpha_1\alpha_2\alpha_3 \in \widetilde{A}^\#$ and similarly for $\widetilde{B}^\#$. By [AnBrRz, Theorem IV.7.2] this means that $\widetilde{A}^\#$ and $\widetilde{B}^\#$ are principal (in $X$), say $\widetilde{A}^\# = \{g' > 0\}$ and $\widetilde{B}^\# = \{g' < 0\}$. Now take any function $g \in K$ whose restriction to $X$ is $g'$. Obviously $g$ separates $\widetilde{A}^\#$ and $\widetilde{B}^\#$.

Finally the implication from *b)* to *a)* is trivial. □

The following result translates the notion of odd and even walls into spaces of orderings. In the proof we use the immediate fact that for any $\widetilde{A} \subset \Sigma_W[\mathbb{Z}_2]$ we have $(\widetilde{A}_W)^\# = (\widetilde{A}^\#)_W$.

**Proposition 4.8.** a) $W$ *is odd if and only if there is* $\sigma \in \Sigma_W[\mathbb{Z}_2]$ *such that* $\sigma \in \widetilde{A}^\#$ *and* $i\sigma \in B^\#$.
  b) $W$ *is even if and only if there is* $\sigma \in \Sigma_W[\mathbb{Z}_2]$ *such that* $\sigma \in \widetilde{A}^\#$ *and* $i\sigma \in \widetilde{A}^\#$ *or* $\sigma \in B^\#$ *and* $i\sigma \in B^\#$.

*Proof.* We only prove *a)*. the proof of *b)* is analogous. So, suppose that $W$ is odd. Replacing $A$ and $B$ by the sets $A'$ and $B'$ stated in the definition of odd walls, we may assume that $A$ and $B$ are separable and that any function separating them changes sign at $W$. Let $f$ be any such function, say with $f|_A > 0$, $f|_B < 0$. In particular $f$ vanishes on $W$. If $\widetilde{A}_W \cap \widetilde{B}_W \neq \emptyset$ then *a)* follows at once. So, suppose that $\widetilde{A}_W \cap \widetilde{B}_W = \emptyset$. Assume that $\widetilde{A}_W$ and $\widetilde{B}_W$ are separable, say $g|_{\widetilde{A}_W} > 0$, $g|_{\widetilde{B}_W} < 0$, where $g \in \Re(M)$ is a unit in $\Re(M)_{\mathfrak{J}(W)}$. The idea is that then $g$ will separate $A$ and $B$ in a small neighbourhood of $W$, contradicting the hypothesis that $W$ is odd.

To make this idea precise, consider the sets $\{g > 0\}$ and $\{g < 0\}$. We claim that the closure of the open semialgebraic set $T = (\{g > 0\} \cap B) \cup (\{g < 0\} \cap A)$ intersects $W$ into a subset of codimension at least one, or in other words that this closure does not contains any open set of the regular locus of $W$. For, otherwise there would be an ordering $\tau \in \widetilde{T}$ which would specialize to an ordering $\gamma \in \Sigma_W$. Now, if $\tau \in \{g > 0\} \cap \widetilde{B}$ we get $\gamma \in \widetilde{B}_W$ and $g(\gamma) > 0$, a contradiction, and similarly if $\tau \in \{g < 0\} \cap \widetilde{A}$. Finally we apply Lojasiewicz's inequality to the functions $f$, $g$ and the set $\overline{T}$, so that there exists a nonnegative function $\varepsilon$ such that $F = f + \varepsilon g$ has the same sign that $f$ on $\overline{T}$ and the zero set of $\varepsilon$ is included in $\overline{T} \cap \{f = 0\}$. In particular, by the remark just made we have that $\varepsilon$ does not vanish identically on $W$. Obviously $F$ separates $A$ and $B$ but it is unit in $\Re(M)_{\mathfrak{J}(W)}$ against our hypothesis on the odd character of $W$.

Coming back to our proof, we may thus assume that $\widetilde{A}_W$ and $\widetilde{B}_W$ are not separable. By Proposition 4.7 this means that there is an ordering $\gamma \in \Sigma_W$ with $\gamma \in (\widetilde{A}_W)^\# \cap (\widetilde{B}_W)^\# = \widetilde{A}_W^\# \cap \widetilde{B}_W^\#$. Since $A$ and $B$ are separable we have $\widetilde{A}^\# \cap \widetilde{B}^\# = \emptyset$, and taking the two generizations of $\gamma$ we get the two orderings $\sigma$ and $i\sigma$ of the statement, so that the proof is complete. □

Finally we combine the above results to prove the main result of this section:

**Theorem 4.9.** (Universal obstruction to separation) *$A$ and $B$ are not generically separable if and only if in some model (obtained from $M$ by a finite sequence of blowings-up) there is a wall of $A$ and $B$ which is simultaneously odd and even.*

*Proof.* If $A$ and $B$ are generically separable it is immediate that there are no walls which are odd and even. So, assume that $A$ and $B$ are not generically separable. Then there exists a finite subspace of orderings $X \subset \Sigma_K$ such that $\widetilde{A} \cap X$ and $\widetilde{B} \cap X$ are not separable. Take $X$ with minimal cardinal. By Theorem 3.1 we may assume that $X$ is a geometric space of orderings. Furthermore, by the minimality assumption we have that $X$ is an extension, say $X = Y[\mathbb{Z}_2]$, with $Y$ a geometric space of orderings. Now, take a model $N$ of $K$ in which $X$ is fully realized, cf. Proposition 2.5. Then $Y$ is a GSO of a hypersurface $W$ of $N$ which by Proposition 4.3 is a wall of $A$ and $B$. We claim that $W$ is an odd and even wall.

Indeed, again by the minimality assumption on $X$ we get that $X \subset \widetilde{A}^\# \cup \widetilde{B}^\#$. Moreover, since $\widetilde{A} \cap X$ and $\widetilde{B} \cap X$ are not separable there is $\sigma \in X$ with $\sigma \in \widetilde{A}^\# \cap \widetilde{B}^\#$. Now $i\sigma \in X$, so that either $i\sigma \in \widetilde{A}^\#$ or $i\sigma \in \widetilde{B}^\#$. In any case, by the above proposition, we get that $W$ is odd and even. □

**Remark 4.10** *a)* Notice that Theorem 4.9 shows that any GSO of minimal cardinal in which $\widetilde{A}$ and $\widetilde{B}$ are not separated is centered at an odd and even wall in a model where it is fully realized.

*b)* Let us derive the following simple fact from the above proof which will be used later: if $A$ and $B$ are not generically separated, there is a fan $F = \{\alpha, i\alpha, \beta, i\beta\} \subset X$ such that, up to renaming, we have: $\alpha \in \widetilde{A}$, $\beta \in \widetilde{B}$ and $\alpha, i\alpha, i\beta \in \widetilde{A}^\#$. Indeed, continuing with the argument at the end of the proof of the Theorem, we have found $\sigma \in \widetilde{A}^\# \cap \widetilde{B}^\#$ with $i\sigma \in \widetilde{A}^\#$ or $i\sigma \in \widetilde{B}^\#$. Suppose the latter. Now, let $\beta \in \widetilde{B} \cap X$ and $\alpha \in \widetilde{A} \cap X$. From $i\beta = \sigma i \sigma \beta$ we get that $i\beta \in \widetilde{A}^{\#}$, and from $i\alpha = \sigma i \sigma \alpha$ we get that $i\alpha \in \widetilde{A}^{\#}$ what proves the claim.

*c)* The following simple example shows that the consideration of blowings-up of $M$ is indeed necessary to find odd an even walls. Let $A$ and $B$ as in the figure:

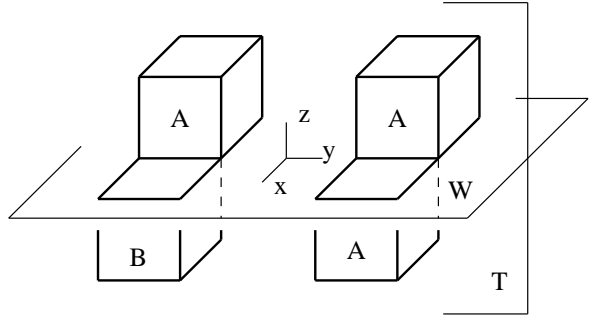

$A$ and $B$ are not separable (in fact they are not separable neither in $\Sigma_W[\mathbb{Z}_2]$ nor in $\Sigma_T[\mathbb{Z}_2]$, where $W$ and $T$ are the planes (walls) shown in the figure). However, neither $W$ nor $T$ are odd or even, as a quick inspection shows. Now consider the quadratic transform given by the equations $z = z'x'$, $y = y'$ and $x = x'$. Then our picture transforms into:

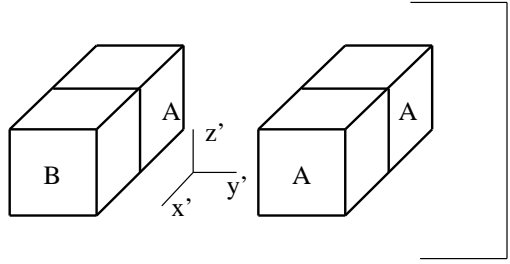

and the wall $\{x' = 0\}$ is obviously odd and even. This example shows also that for a wall $W$ the statements "$A$ and $B$ are not separable in $\Sigma_W[\mathbb{Z}_2]$" and "$W$ is odd and even", are not equivalent. Obviously the latter implies the former, but the converse is not true in general: in the first picture above, $A$ and $B$ are not separable in $\Sigma_W[\mathbb{Z}_2]$, but $W$ is neither odd nor even.

*d)* Notice that if for a wall $W$ we have $\text{Reg}W \cap (A \cup B) \neq \emptyset$ then $W$ is even (remember that $A$ and $B$ are open). Indeed if $a$ is a regular point of $W$ which lies, say, in $A$ then a whole neighbourhood of $a$ is contained in $A$, which obviously implies that $W$ is even. □

## 5. The geometric criterion

Despite the beauty of the universal obstruction proved above, in order to decide the generic separation of two semialgebraic sets we want to fix a model of $M$ and develop a test in it, rather that consider all walls in all possible models of $K$. It turns out that the key property for $A$ and $B$ which we need to do it (besides the non-singularity of $M$) is that their walls be normal crossings. Since we will work algebraically, by this we mean that *if $f_1, \ldots, f_s$ is a family of polynomials describing*

*A and B, then all the $f_i$'s are normal crossings.* We start by showing that under this condition we do not, any longer, need to leave $M$:

**Theorem 5.1.** *Assume that $M$ is non singular and that all the walls of $A$ and $B$ are nonsingular and normal crossings. Then $A$ are $B$ are generically separable if and only if $\widetilde{A}$ and $\widetilde{B}$ are separable in all subspaces $\Sigma_W[\mathbb{Z}_2]$ for all the walls $W$ of $A$ and $B$ in $M$.*

*Proof.* Suppose that $A$ and $B$ are not generically separable. We may find a GSO, say $X$, of minimal cardinal in which $\widetilde{A}$ and $\widetilde{B}$ are not separable. Then by Proposition 2.5, Theorem 4.9 and Remark 4.10 *a)* there is a model $M'$ of $K$, obtained blowing-up $M$, say $\pi: M' \to M$, such that $X$ is centered at a wall $W'$ of $A' = \pi^{-1}(A)$ and $B' = \pi^{-1}(B)$ in $M'$ which is odd and even. This means that there is a 4-element fan $F = \{\sigma, i\sigma, \tau, i\tau\}$ in $\Sigma_{W'}[\mathbb{Z}_2] \subset \Sigma_K$ with, say, $\sigma, i\sigma, i\tau \in \widetilde{A}^\#$ and $\tau \in \widetilde{B}$ and $\sigma \in \widetilde{A}$, cf. Remark 4.10 *b)*. Set $R' = \mathfrak{R}(M')$, $\mathfrak{q} = \mathfrak{J}(W')R'$, $R = \mathfrak{R}(M)$, $T = \pi(W')$ and $\mathfrak{p} = \mathfrak{q} \cap R$, so that $\mathfrak{p} = \mathfrak{J}(T)R$. We have that $R'_\mathfrak{q}$ is a rank 1 discrete valuation ring with residue field $K(W')$ and $R_\mathfrak{p}$ is a regular local ring of dimension, say $d$. If $d = 1$ then $\Sigma_{W'} = \Sigma_T$, $T$ is again an odd and even wall and we are done. So, suppose that $d \geq 2$. Let $f_1, \ldots, f_s$ be a family of polynomials describing $A$ and $B$. Since all the $f_i$'s are normal crossings, there is a regular system of parameters $x_1, \ldots, x_d$ of the ring $R_\mathfrak{p}$ such that for all $i = 1, \ldots, s$ we have
$$f_j = a_j x_1^{m_1} \cdots x_d^{m_d}.$$
where $a_j$ is a unit in $R_\mathfrak{p}$.

Changing $x_i$ by $-x_i$ we may assume that $x_i(\sigma) > 0$ for all $i$. Now, let $t$ be a uniformizing parameter of the valuation ring $R'_\mathfrak{q}$, so that, say $t(\sigma) > 0$ and $t(\tau) > 0$. For each $i = 1, \ldots, d$ we have $x_i = t^{r_i} v_i$ where $v_i$ is a unit in $R'_\mathfrak{q}$. We claim that some of the exponents $r_i$ is odd. For otherwise, for all $j$ and all $\sigma \in \Sigma_{W'}[\mathbb{Z}_2]$, the function $f_j$ would have the same sign in $\sigma$ and $i\sigma$. This would imply that $\widetilde{A}$ and $\widetilde{B}$ are saturated with respect to the specialization map $\Sigma_{W'}[\mathbb{Z}_2] \to \Sigma_{W'}$ which in turn contradicts the minimality of $X$. Hence, assume that $r_1$ is odd, so that $x_1$ separates $\sigma$ from $i\sigma$. To keep the geometrical intuition, let us represent by $W$ the wall defined by $\{x_1 = 0\}$ in $M$.

Now we distinguish two cases: assume first that the specializations $\overline{\sigma}$ and $\overline{\tau}$ of $\sigma$ and $\tau$ in $W'$ are separated by some $a_i$, say $a_1$ (what in particular implies that they are already different in $T$). Then, to any $\gamma \in (\widetilde{A} \cup \widetilde{B}) \cap \Sigma_{W'}[\mathbb{Z}_2]$ we attach an ordering in $\Sigma_W[\mathbb{Z}_2]$ constructed as follows: $R_{(x_1)}$ is a discrete valuation ring with residue field the quotient field of $R/(x_1)$. In turn, $(R/(x_1))_\mathfrak{p} = R_\mathfrak{p}/(x_1)$ is a regular local ring of dimension $d - 1$ with $x_2, \ldots, x_d$ as parameters. Thus, we may consider the rank $d - 1$ discrete valuation $v_{d-1}$ of the quotient field of $(R/(x_1))$ introduced in the Construction 1.1 associated to this system of parameters and dominating $R_\mathfrak{p}/(x_1)$. Let $v_d$ be the composite of $R_{(x_1)}$ and $v_{d-1}$. In particular $v_d$ is a rank $d$ discrete valuation ring of $K$ with residue field $\kappa(\mathfrak{p})$, and for any ordering of $\Sigma_{\kappa(\mathfrak{p})}$ their pullbacks by $v_d$ are completely determined by the signs of $x_1 \ldots, x_d$ and belong to $\Sigma_W[\mathbb{Z}_2]$.

Thus, given any ordering $\gamma \in \Sigma_{W'}[\mathbb{Z}_2]$ we define $\widehat{\gamma}$ as the pull-back of $\overline{\gamma} \cap \kappa(\mathfrak{p})$ by $v_d$ defined by $x_i(\widehat{\gamma}) = x_i(\gamma)$. Notice that the correspondence $\gamma \to \widehat{\gamma}$ is not in general bijective. However, for any $j = 1, \ldots, s$ we have $f_j(\gamma) = f_j(\widehat{\gamma})$, so that if $\gamma$ was in $\widetilde{A}$ (resp. $\widetilde{B}$), so is $\widehat{\gamma}$. Now, an immediate checking shows that $\widehat{\sigma}, \widehat{i\sigma}, \widehat{\tau}, \widehat{i\tau}$ is a fan of $\Sigma_W[\mathbb{Z}_2]$ with $\widehat{\sigma} \in \widetilde{A}$, $\widehat{\tau} \in \widetilde{B}$ and $\widehat{\sigma}, i\widehat{\sigma}, i\widehat{\tau} \in \widetilde{A}^{\#}$. This shows that $\widetilde{A}$ and $\widetilde{B}$ are not separable in $\Sigma_W[\mathbb{Z}_2]$ for the wall $W$, against our assumption.

Secondly, suppose now that $\overline{\sigma}$ and $\overline{\tau}$ cannot be separated by any $a_i$ (for instance if $\overline{\sigma} = \overline{\tau}$). Then, since $\sigma$ and $\tau$ are separated by some $f_j$, there is another $x_i$, say $x_2$ such that $x_2(\tau) < 0$. Thus, we have the following distribution of signs for $x_1$ and $x_2$:

|       | $\sigma$ | $i\sigma$     | $\tau$ | $i\tau$        |
|-------|----------|---------------|--------|----------------|
| $x_1$ | 1        | $-1$          | 1      | $-1$           |
| $x_2$ | 1        | $\varepsilon$ | $-1$   | $-\varepsilon$ |

where $\varepsilon = 1$ or $-1$. In any case we discover the shape of a four element fan, and we proceed as above: for any $\gamma \in \Sigma_{W'}[\mathbb{Z}_2]$ we define $\widehat{\gamma}$ as the pull-back of $\overline{\gamma} \cap \kappa(\mathfrak{p})$ by $v_d$ defined by $x_i(\widehat{\gamma}) = x_i(\gamma)$. As before, an immediate checking shows that $\widehat{\sigma}, i\widehat{\sigma}, \widehat{\tau}, i\widehat{\tau}$ is a fan of $\Sigma_W[\mathbb{Z}_2]$ with $\widehat{\sigma} \in \widetilde{A}$, $\widehat{\tau} \in \widetilde{B}$ and $\widehat{\sigma}, i\widehat{\sigma}, i\widehat{\tau} \in \widetilde{A}^{\#}$. Therefore we get that $A$ and $B$ are not separable in $\Sigma_W[\mathbb{Z}_2]$ for the wall $W$, against our assumption, and the proof is complete. $\square$

Once we know that we have to check only the separation locally at the walls in a fixed model, the obvious question is, how can we test if $A$ and $B$ are separable in $\Sigma_W[\mathbb{Z}_2]$? Fortunately there is a very nice geometric criterion to do that in terms of the shadows and counter–shadows of $A$ and $B$ in $W$. We define the *shadows* of $A$ and $B$ in $W$ as the sets $A_W = \overline{A} \cap W$ and $B_W = \overline{B} \cap W$. We define the *counter–shadows* of $A$ and $B$ in $W$ as the shadows of the sets

$$A^t = (A \cap \{t > 0\}) \cup (B \cap \{t < 0\})$$
$$B^t = (A \cap \{t < 0\}) \cup (B \cap \{t > 0\})$$

where $t$ is a uniformizer of the valuation ring $\mathfrak{R}(M)_{\mathfrak{J}(W)}$. Notice that changing $t$ by $-t$ reverses the roles of $A^t$ and $B^t$. Thus, since we are working generically, up to reversing the roles of $A$ and $B$, the counter–shadows do not depend on the uniformizer.

The following immediate remark shows some useful relations:

**Proposition 5.2.** *a) $f$ separates generically $A$ and $B$ if and only if $tf$ separates generically $A^t$ and $B^t$.*
*b) $W$ is an even wall for $A$ and $B$ if and only if it is odd for $A^t$ and $B^t$.* $\square$

The idea is that if the shadows can be generically separated, then $\widetilde{A}$ and $\widetilde{B}$ can be separated in $\Sigma_W[\mathbb{Z}_2]$ and therefore $A$ and $B$ can be generically separated in a (infinitesimal) neighbourhood of $W$. However the converse is not true: take $A = \{z > 0, y + 1 > 0\}$ and $B = \{z < 0, 1 - y > 0\}$. Their shadows in $z = 0$ cannot be separated, but still the function $z$ separates $A$ from $B$. Here is where the counter-shadows give a hand: they take care of this phenomena, and as Theorem below shows, shadows and counter–shadows together characterize the separability in $\Sigma_W[\mathbb{Z}_2]$.

The shadow of a semialgebraic set $A$ corresponds in the context of a geometric spaces of orderings $X = X'[\mathbb{Z}_2]$ to the set of specializations of $\widetilde{A}$ in $X'$. Similarly, given two constructible subsets $\widetilde{A}, \widetilde{B} \subset X$, the counter–shadows are the specializations in $X'$ of the sets:
$$\widetilde{A}^t = (\widetilde{A} \cap (X' \times \{1\})) \cup (\widetilde{B} \cap (X' \times \{i\}))$$
$$\widetilde{B}^t = (\widetilde{A} \cap (X' \times \{i\})) \cup (\widetilde{B} \cap (X' \times \{1\}))$$

We have the following general result:

**Proposition 5.3.** $\widetilde{A}$ and $\widetilde{B}$ are separable in $X = X'[\mathbb{Z}_2]$ if and only if either their shadows or their counter–shadows are separable in $X'$.

*Proof.* If the shadows are separable, say by a function $g'$ on $X'$, then the same function separates $\widetilde{A}$ and $\widetilde{B}$. Analogously, if the counter–shadows are separable then $\widetilde{A}^t$ and $\widetilde{B}^t$ are separable and therefore so are $\widetilde{A}$ and $\widetilde{B}$ by Proposition 5.2. Conversely, assume that $g$ separates $\widetilde{A}$ and $\widetilde{B}$ and denote by $\widetilde{A}'$ and $\widetilde{B}'$ respectively their specializations in $X'$. If $g = (1, g')$ then $g'$ separates $\widetilde{A}'$ and $\widetilde{B}'$. If $g = (i, g')$, then $ig = (1, g')$ separates $\widetilde{A}^t$ and $\widetilde{B}^t$ and therefore $g'$ separates the counter–shadows of $\widetilde{A}$ and $\widetilde{B}$ in $X'$. □

As a consequence we get the following

**Corollary 5.4.** Let $A$ and $B$ be open semialgebraic. $\widetilde{A}$ and $\widetilde{B}$ are not separable in $\Sigma_W[\mathbb{Z}_2]$ if and only if neither the shadows nor their counter–shadows of $A$ and $B$ in $W$ are separable in $\Sigma_W$. □

Winding up, taking into account Theorem 5.1 above we get the following nice geometric criterion for separation:

**Theorem 5.5.** (Separation criterion) *Assume that $M$ is non singular and $A$ and $B$ have walls at normal crossings. Then $A$ and $B$ are generically separable if and only if for any wall $W$ of $A$ and $B$ in $M$, either the shadows or the counter–shadows of $A$ and $B$ are separable in $W$.*

*Proof.* By Theorem 5.1 $A$ and $B$ are separable if and only if they are separable in the spaces $\Sigma_W[\mathbb{Z}_2]$ for all walls $W$, and Corollary 5.4 states that this is equivalent to say that either the shadows or the counter–shadows of $A$ and $B$ in $W$ are separable. So we are done. □

**Remark 5.6** *a)* The condition of the walls of $A$ and $B$ be at normal crossings is important. Notice that in the example in the introduction the shadows and counter–shadows of $A$ and $B$ are separated when they intersect only at the origin, while they are not after blowing-up.

*b)* The geometric criterion translates the problem of generic separation of $A$ and $B$ into a finite number of similar problems in one less dimension. This provides a recursive method to check generic separation, which in the end makes it a decidable problem by using Tarski's Principle as we will see in the next section. □

## 6. Decision procedure for generic separation

We want to prove that the generic separation of two open semialgebraic sets $A$ and $B$ on an algebraic variety $M \subset \mathbb{R}^n$ is decidable. Since we are concerned here with effective algorithms, we assume that $M$, as well as $A$ and $B$, are defined over $\mathbb{R}_0$, the field of real algebraic numbers. It is clear from Tarski's Transfer Principle that all the semialgebraic sets constructed in the sequel (like the irreducible components of the Zariski closure of the boundary of $A$ and $B$) can be defined as well over $\mathbb{R}_0$. Under this hypothesis we have:

**Theorem 6.1.** *The generic separation of semialgebraic sets on an affine variety over the real numbers is decidable.*

We will present two proofs. The first presents a tentative algorithm to test the separability, having $M$, $A$ and $B$ as input and answering YES or NO. It relies on the fact announced in [BiMi] that the desingularization can be perfomed algorithmically. The rest of the operations needed (as computing closures of semialgebraic sets, computing the dimension of an algebraic set or finding its irreducible components) belong to the yoga of semialgebraic geometry and algebraic geometry algorithms, see [He], [BeNe], [Ne]. The second proof is more in the spirit of Model Theory and uses the notion of recursive enumerability. Both proofs work by induction on $d = \dim M$.

Notice that for $d = 1$, $A$ and $B$ are generically separable (and separable since both notions coincide) if and only if $A \cap B = \emptyset$, which is a statement that can be expressed by a (formal) sentence in the first order language of ordered fields and therefore its truth in $\mathbb{R}$ is decidable by Tarski's result. Thus we assume that $\dim M = d > 1$.

*First proof:*

**Generic separation Algorithm.** Having $M$, $A$ and $B$ as input answers YES or NO to the question of whether they are generically separable:

1. Find a desingularization $\pi : M' \to M$ and put $\overline{\partial A}^z \cup \overline{\partial B}^z$ into normal crossings. Set $M \leftarrow M'$, $A \leftarrow \pi^{-1}(A)$, $B \leftarrow \pi^{-1}(B)$.
2. For all walls $W_1, \ldots, W_r$ of $A$ and $B$ check whether the shadows $A_{W_i}$, $B_{W_i}$ and the counter–shadows $A_{W_i}^t$, $B_{W_i}^t$ are generically separable (apply induction). If for some $i$ both are not generically separable return "NO, $A$ and $B$ are not generically separable". Otherwise return "YES, they are generically separable".

*Second proof:* Again we work by induction on $\dim M$, so that we assume that we have proved decidability for generic separation in all varieties of dimension less than $d$.

Decidability of generic separation in $M$ with $\dim M = d$ will now be proved by recursively listing, on the one hand, all triples $(M, A, B)$ where $M$ is a variety of dimension $d$ and $A$ and $B$ are open semialgebraic subsets of $M$ which can be generically separated in $M$, and on the other hand, we recursively list those triples where $A$ and $B$ cannot be generically separated in $M$.

Thus given a particular variety $M_0$ and a pair $A_0$, $B_0$ of open semialgebraic sets in $M_0$ the decision about (generic) separation is obtained by passing simultaneously through both lists. Depending on which side the triple $(M_0, A_0, B_0)$ will occur we know whether $A_0$ and $B_0$ are generically separable in $M_0$ or not.

The first list for generic separation is clearly contained in the list of quadruples $(M, A, B, p)$ of varieties $M$ of dimension $d$ and open semialgebraic sets $A$ and $B$ in $M$ which are generically separated by the polynomial $p$. This list of quadruples may be obtained in the following way. First, notice that by Bröcker's results, cf. [Br4], [AnBrRz], the number of polynomials needed to describe any semialgebraic subset of $M$ is bounded by a constant depending only on the dimension $d$ of $M$, so that we may recursively list all quadruples $(M, A, B, p)$ of all semialgebraic open sets $A$ and $B$ and all polynomials $p$. Then we use Tarski's decision procedure to eliminate those quadruples for which $p$ does not generically separate $A$ from $B$.

The second list will be contained in the following recursive list of certain finite sequences:
$$(M, A, B, X, (\pi_\nu)_{\nu \leq m}, M', A', B', (W'_j, g_j)_{j \leq s}, X')$$
where we have:
- $M \subset \mathbb{R}^n$ and $M' \subset \mathbb{R}^m$ are compact varieties of dimension $d$,
- $X$ and $X'$ are Zariski closed proper subsets of $M$ and $M'$ respectively,
- $\pi = (\pi_\nu)_{\nu \leq n}$ is a biregular isomorphism of $M' \setminus X'$ and $M \setminus X$,
- $A$ and $B$ are open semialgebraic sets in $M$,
- $A'$ and $B'$ are open semialgebraic sets in $M'$,
- $W'_1, \ldots, W'_s$ are the irreducible components of $\overline{\partial A'}^Z \cup \overline{\partial B'}^Z$,
- $g_1, \ldots, g_s$ are polynomials changing sign on $W'_1, \ldots, W'_s$ respectively.

Using the fact that it is decidable whether finitely many given algebraic sets form the components of the Zariski closure of a given semialgebraic set, cf. [He], [BeNe], [Ne], it is clear that such a list can be recursively obtained.

Now we shorten this list so that eventually only the triples $(M, A, B)$ where $A$ and $B$ are not generically separable remain. We do this first by eliminating all sequences not satisfying one of the following conditions, which are, all of them, decidable.

We require that:
(i) $M'$ is nonsingular
(ii) $A' = Int(\overline{\pi^{-1}(A \setminus X)})$, $B' = Int(\overline{\pi^{-1}(B \setminus X)})$
(iii) $W'_1, \ldots, W'_s$ are nonsingular and normal crossings.

Finally we save only those sequences which do not satisfy the separation criterion, that is, which satisfy the following condition which is decidable by the induction hypothesis, so that the procedure for generic separtion is complete:

(a) for some $j \leq s$, the shadows of $A'$ and $B'$ as well as the shadows of $(A')^{g_j}$ and $(B')^{g_j}$ on $W'_i$ are not generically separable, □

We may ask whether given two generically separable semialgebraic sets $A$ and $B$ there is an upper bound on the degree of a separating polynomial in terms of their "complexity" that is, the number of polynomials, their total degree and size of their coefficients. First of all, by taking limits, one can see that if $A$ is an open semialgebraic set and $\{B_n\}$ is a properly nested family, $B_n \subset B_{n+1}$, of open semialgebraic sets such that for all $n$, $A$ and $B_n$ are separable while $A$ and $B = \bigcup_n B_n$ are not generically separable, then the degrees of any sequence $\{f_n\}$ of polynomials separating $A$ and $B_n$ must go to infinity. Indeed, notice that scaling by the coefficient with larger absolute value we may always assume that the coeficients of $f_n$ lie in the unit ball. Then, if the degrees were bounded, say by $\rho \in \mathbb{N}$, we might represent each $f_n$ as a $\rho$–tuple in $[-1,1]^\rho$ and we would get a limit point $f$ which would represent a polynomial separating $A$ and $B$, a contradiction.

Also, by model completeness, it follows that it cannot exist a recursive bound for the degree of the separating function, since then the separation criterion 5.5 would be true for any real closed field $R$, while the following example shows that it fails when $R$ is non-archimedean.

**Example 6.2** Let $R$ be any non-archimedean real closed field, and consider the open semialgebraic subsets $A, B \subset R \times R$ as in the left hand side of the following picture, where $\varepsilon \in R$ is infinitesimal, i.e. $\varepsilon < \dfrac{1}{n}$ for all $n \in \mathbb{N}$.

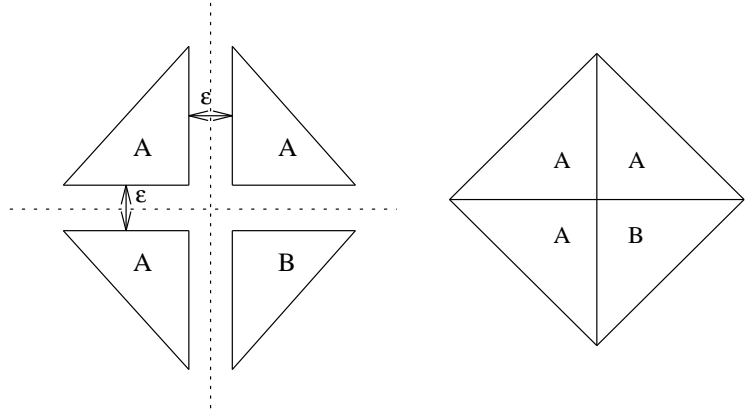

Notice that $A$ and $B$ verify the hypothesis of Theorem 5.5, since $R \times R$ is nonsingular and the walls are nonsingular and normal crossings. However we claim that there is no polynomial separating $A$ and $B$, even generically.

For suppose that $h \in R[x,y]$ separates generically $A$ and $B$, and let $\mathfrak{O}$ be the convex hull of $\mathbb{Q}$ in $R$. then $\mathfrak{O}$ is a convex valuation of $R$ whose maximal ideal $\mathfrak{m}$ is the set of infinitesimal elements. In particular $\varepsilon \in \mathfrak{m}$ and the residue field $R'$ of $\mathfrak{O}$ is an archimedean real closed field, so that it can be considered as a subfield of $\mathbb{R}$. Moreover, by scaling $h$ with the coefficient of larger absolute value in $R$ we may assume that $h \in \mathfrak{O}[x,y]$ and not all coefficients of $h$ lie in $\mathfrak{m}$. Now, the image $A'$, $B'$ of $A$ and $B$ in $R'$ under the residue map looks like the right hand side of the previous picture, and it is obvious that they are not generically separable in $R'$. However, if $h'$ is the image of $h$ in $R'$ we have $h' \neq 0$ and that it is, say $\geq 0$ on $A'$ and $\leq 0$ on $B'$, so that we get $h' \in \mathfrak{m}[x,y]$, a contradiction. □

Since besides Tarski's Principle the proof of the criterion depends on the density theorem of geometric spaces of orderings, we see that this approximation result holds only over $\mathbb{R}$. In fact, in the example above, there is a 4-element fan $F$ of $\Sigma_{\mathbb{R}^2}$ in which the separation of $\widetilde{A}$ and $\widetilde{B}$ fails, which cannot be approximated by any geometrical one.

## 7. Generic separation versus separation

Having settled the question of whether two semialgebraic sets are generically separable we want to study now when they can be separated. So, in this section $A$ and $B$ are arbitrary semialgebraic subsets of and algebraic variety $M$. We recall that the boundary $\partial A$ of $A$ is defined as $\partial A = \overline{A} \setminus \text{Int}(A)$. In particular, all points where $\dim A < \dim M$ belong to $\partial A$.

Obviously, separation implies generic separation, but the converse is not true as it is easily seen by taking for instance a cubic with an isolated point:

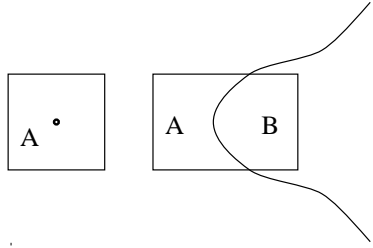

or the sets in $\mathbb{R}^2$

$$A = \{(x,y) \mid 0 < x < 1, 0 < y < 1\} \cup \{(x,y) \mid 0 < x < \frac{1}{2}, y = 0\}$$
$$B = \{(x,y) \mid 0 < x < 1, -1 < y < 0\} \cup \{(x,y) \mid \frac{1}{2} < x < 1, y = 0\}$$

In both cases the idea is the same: any function separating generically $A$ and $B$ must vanish at certain points which are in $A \cup B$ and therefore they cannot be separated in the sense of the definition 4.1. Thus we want to identify the set of common zeros of the family of functions which separate generically $A$ and $B$.

To start with, notice that any function $f$ which separates $A$ and $B$ generically must vanish identically on $\overline{A} \cap \overline{B}$ and therefore in $\overline{\overline{A} \cap \overline{B}}^z$. In fact, this is only a part of the set of points in which any function separating $A$ and $B$ must vanish. For, let $W$ be any irreducible algebraic subset of $M$, such that the shadows of $A$ and $B$ in $W$ are not generically separable. Then any function $f$ which separates $A$ and $B$ must vanish identically on $W$ since it must be $\geq 0$ on $A_W$ and $\leq 0$ on $B_W$. We denote by $Z$ the union of all these sets and we name it *the separation nullspace of $A$ and $B$*. Notice that in particular we have $\overline{\overline{A} \cap \overline{B}}^z \subset Z$ and that a necessary condition for $A$ and $B$ to be separable is that $Z \cap (A \cup B) = \emptyset$.

The following result extends Theorem 3.3 of [AcBgFo] and shows that separation and generic separation coincide up to the separation nullspace $Z$. Compare also with [AnBrRz, Theorem V.3.2 a), b)].

**Theorem 7.1.** *With the notations above, assume that $Z \cap (A \cup B) = \emptyset$. Then $A$ and $B$ are separable if and only if they are generically separable.*

*Proof.* One direction is trivial. For the other, we assume that there exists $f$ which separates generically $A$ and $B$ in $M$, that is, there is a proper algebraic set $W$ such that $f|_{A \setminus W} > 0$ and $f|_{B \setminus W} < 0$. We will see that we can find $f'$ which separates $A$ and $B$ up to a proper algebraic subset $W' \subset W$. Thus the proof ends by induction.

First of all notice that from our assumption on $Z$ we may assume that no irreducible component of $W$ lies in $Z$. Also, replacing $W$ by $\overline{\{f = 0\} \cap (A \cup B)}^z$, we may assume that $W = \overline{W \cap (A \cup B)}^z$ and that $f$ vanishes on $W$. Now, take any irreducible component of $W$ of maximum dimension, which we still denote by $W$. Since $W$ is not in $Z$, there is a function $g$ which separates generically the shadows $A_W$ and $B_W$, say $g|_{(A_W \setminus W')} > 0$ and $g|_{(B_W \setminus W')} < 0$. Consider the semialgebraic set $S = (\overline{A} \setminus \{g > 0\}) \cup (\overline{B} \setminus \{g < 0\})$. By Hörmander-Łojasiewicz inequality there is a nonnegative function $\varepsilon$ such that $f' = f + \varepsilon g$ has the same sign as $f$ over $S$, and with the zero set of $\varepsilon$ contained in $\overline{\{f = 0\} \cap S}^Z$. Thus we get that $f'|_{A \setminus W'} > 0$ and $f'_{|B \setminus W'} < 0$ and we are done. $\square$

**Remark 7.2** *a)* If $\dim M = 2$ then $Z = \overline{\overline{A} \cap \overline{B}}^Z$. If $\dim M = 3$ and $M$ is nonsingular, in [AcBrFo] it is shown that $Z$ is the union of $\overline{\overline{A} \cap \overline{B}}^Z$ and the *odd walls* of $A$ and $B$. However in higher dimensions this is no longer true as the following example shows: let

$$\begin{aligned} A &= \{-1 < x < 0, -1 < y < 0, -1 < z < 1, -1 < t < 0\} \cup \\ &\quad \{-1 < x < 0, -1 < y < 0, -1 < z < 0, 0 < t < 1\} \\ B &= \{0 < x < 1, 0 < y < 1, 0 < z < 1, 0 < t < 1\} \end{aligned}$$

Notice that the shadows of $A$ and $B$ in $W = \{x = y = 0\}$ are the sets of the picture below, which are not generically separable. Thus, $Z$ contains $W$, although none of the walls of $A$ and $B$ are odd, as we will see later.

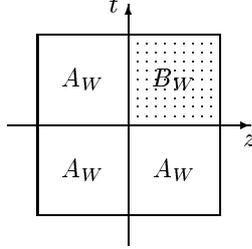

*b)* Theorem 7.1 shows that if we consider the sets $A^* = \overline{A} \cup Z$ and $B^* = \overline{B} \cup Z$, then $A$ and $B$ are separable if and only if $A^*$ and $B^*$ are separable in the sense of [AnBrRz, Definition V.3.1].

*c)* Also Theorem 7.1 shows that $Z$ coincides with the set of common zeros of all functions separating generically $A$ and $B$. This is the reason why we call $Z$ the *separation nullspace* of $A$ and $B$ in $M$. $\square$

The next result shows that $Z$ is contained in the union of $\mathrm{Sing} M$ and the Zariski closure of the boundary of $A \cup B$.

**Theorem 7.3.** *Let $A$ and $B$ be two semialgebraic subsets of $M$ such that $(\overline{\partial A}^Z \cup \overline{\partial B}^Z \cup \mathrm{Sing} M) \cap (A \cup B) = \emptyset$. Then $A$ and $B$ are separable if and only if they are generically separable. In particular $Z \subset \overline{\partial A}^Z \cup \overline{\partial B}^Z \cup \mathrm{Sing} M$.*

*Proof.* Suppose that $A$ and $B$ are not separable but they are generically separable. Then, by Theorem 5.3 there is a subvariety $T \subset Z$ such that $(A \cup B) \cap T \neq \emptyset$ and the shadows $A_T$ and $B_T$ of $A$ and $B$ in $T$ are not generically separable. By Brocker's theorem, this means that there is a finite subspace of orderings $Y \subset \Sigma_T$ such that $\widetilde{(A_T)} \cap Y$ and $\widetilde{(B_T)} \cap Y$ are not separable. In particular we have $A_T \neq \emptyset$ and $B_T \neq \emptyset$. Now, by our hypothesis we get $T \not\subset \overline{\partial A}^Z \cup \overline{\partial B}^Z$ and therefore $A_T \subset A$ and $B_T \subset B$. Let $\mathfrak{p}$ be the ideal of functions vanishing on $T$. Since $T \not\subset \mathrm{Sing} M$, the localization $\mathfrak{R}(M)_{\mathfrak{p}}$ is a regular local ring with residue field $\mathfrak{K}(T)$. Thus, there is a discrete valuation ring $V$ of $K$ dominating $\mathfrak{R}(M)_{\mathfrak{p}}$ and with residue field $K(T)$. Fix any sign for the regular system of parameters of $\Re(M)_{\mathfrak{p}}$ and consider the pull back of $Y$ to $K$ with these parameters fixed. We get a subspace of orderings $X \subset \Sigma_K$ isomorphich to $Y$. Moreover, if $\sigma \in X$ is the unique generizarion of $\tau \in Y$ we have that $\sigma \in \widetilde{A}$ (resp. $\sigma \in \widetilde{B}$) if and only if $\tau \in \widetilde{(A_T)} \cap Y$ (resp. $\tau \in \widetilde{(B_T)} \cap Y$) and we get that $\widetilde{A}$ and $\widetilde{B}$ are not separable in $X$, and therefore neither are they in $\Sigma_K$, in contradiction with our assumption that $A$ and $B$ were generically separable. $\square$

As an application we see that in the normal crossing situation of the separation criterion 5.5 the notions of separation and generic separation coincide for open sets.

**Corollary 7.4.** *Assume that $A$ and $B$ are open, $M$ is non-singular and all the walls are nonsingular and normal crossings. Then $A$ are $B$ are separable if and only if they are generically separable.*

*Proof.* Assume that $A$ and $B$ are generically separable. We will see that then $Z \cap (A \cup B) = \emptyset$, what by Theorem 7.1 implies that $A$ and $B$ are separable. We work by induction on $\dim M$. If $\dim M = 1$ the claim is trivial. So assume that $d = \dim M > 1$ and that $Z \cap (A \cup B) \neq \emptyset$. Let $T \subset Z$ be an irreducible algebraic set with $T \cap (A \cup B) \neq \emptyset$. Since $M$ is nonsingular, by Theorem 7.3 we have that $T \subset \overline{\partial A}^z \cup \overline{\partial B}^z$, so that $T$ is contained in a wall, say $W$. In particular $W \cap (A \cup B) \neq \emptyset$ and since $W$ is nonsingular, by Remark 4.10 d), we have that $W$ is even. Thus, for any function $f$ separating generically $A$ and $B$, in $\Re(M)_{\mathfrak{J}(W)}$ we have $f = t^m u$ where $t$ is a uniformizer, $u$ a unit and $m$ is even. Hence $u$ still separates $A$ and $B$ generically and since it is a unit, it separates also $A_W$ and $B_W$.

Since $W$ is nonsingular and the walls of $A_W$ and $B_W$ are also nonsingular and normal crossings, we may apply the induction hypothesis to $W$, $A_W$ and $B_W$. In particular, since $A_T = (A_W)_T$ and $B_T = (B_W)_T$ we have that $T$ is also in the separation nullspace of $A_W$ and $B_W$ and by induction we get that $T \cap (A_W \cup B_W) = \emptyset$. But $T \cap (A \cup B) \subset T \cap (A_W \cup B_W)$ and we get a contradiction. $\square$

The following example shows that the result of Theorem 7.3 cannot be improved.

**Example 7.5** Let $Q$ be an irreducible compact curve in $\mathbb{R}^2$ with two connected components $Q_1$ and $Q_2$ which are separated by the line $x_1 = 0$. Then consider the cone $M \subset \mathbb{R}^5$ constructed over $Q$ placed at the plane $x_3 = x_4 = x_5 = 1$ and taking the 2-plane with coordinates $x_4, x_5$ as vertex. Thus, $M$ is the union of two closed semialgebraic sets $M = M_1 \cup M_2$, where $M_i$ is the cone over $Q_i$, and the polynomial $x_1 x_5$ separates them outside the vertex $M_0$ of $M$, so that they are generically separated. On the vertex we consider the semialgebraic sets

$$\begin{aligned} A_0 &= \{-1 < x_4 < 0, -1 < x_5 < 1\} \cup \{0 < x_5 < 1, -1 < x_4 < 1\} \\ B_0 &= \{0 < x_4 < 1, -1 < x_5 < 0\} \end{aligned}$$

Next, we define $A$ as the union of all the generatrices of the cone connecting $Q_1$ and $A_0$, and $B$ as the union all the generatrices of the cone connecting $Q_2$ and $B_0$. Then $A$ and $B$ are open semialgebraic sets in $M$, generically separated by $x_1 x_5$ but $A \cap M_0$ and $B \cap M_0$ are not. So, $M_0 = \text{Sing} M = Z$. $\square$

Although Theorem 7.3 does not identify the separation nullspace, once we know that $A$ and $B$ are generically separated, it opens the door to decide whether indeed $Z \cap (A \cup B) = \emptyset$. Indeed, let $T \subset Z$ be an irreducible subvariety with $T \cap (A \cup B) \neq \emptyset$ and maximal dimension. Then $T$ is contained in some irreducible component $W$ of $(\overline{\partial A}^z \cup \overline{\partial B}^z \cup \mathrm{Sing}M)$. In particular $W \cap (A \cup B) \neq \emptyset$. Now, if the shadows $A_W$, $B_W$ are not generically separated, then $W$ itself is in $Z$, $W = T$ and we are done with the test. So assume that $A_W$ and $B_W$ are generically separated. Notice that $(A_W)_T = A_T$, $(B_W)_T = B_T$ and $T \cap (A \cup B) \subset T \cap (A_W \cup B_W)$ since $A \cap W \subset A_W$ and $B \cap W \subset B_W$. Thus in particular $T$ is contained in the separation nullspace of $A_W$ and $B_W$ in $W$ which by Theorem 5.7 is contained in $(\overline{\partial A_W}^z \cup \overline{\partial B_W}^z \cup \mathrm{Sing}W)$. Hence $T$ must be contained in some irreducible component of this set, which has codimension at least 1 in $W$. Iterating this process, we get that $T$ must be one of the irreducible varieties appearing in the finite list $\mathfrak{L}_{M,A,B}$ obtained in this way, which we systematize in the following algorithm:

**Producing the list $\mathfrak{L}_{M,A,B}$ of obstructions to separation:**
1. Initialize $\mathfrak{L}_{M,A,B} \leftarrow \{\emptyset\}$, $M_0 \leftarrow M$, $A_0 \leftarrow A$ and $B_0 \leftarrow B$.
2. Add to $\mathfrak{L}_{M,A,B}$ all the irreducible components $W$ of $\mathrm{Sing}(M_0) \cup (\overline{\partial A_0}^z \cup \overline{\partial B_0}^z)$ which intersect $A \cup B$.
3. For the new $W \in \mathfrak{L}_{M,A,B}$ IF $\dim W > 1$ set $M_0 \leftarrow W$, $A_0 \leftarrow A_W$ and $B_0 \leftarrow B_W$ and GOTO 2.

Finally, bringing together the generic separation decision procedure and the procedure for the list $\mathfrak{L}_{M,A,B}$ we can design a procedure to decide the separation of $A$ and $B$, proving the desired result:

**Theorem 7.6.** *The separation of semialgebraic sets on an affine variety over the real numbers is decidable.*

*Proof.* (**Separation Decision Procedure**) Having $M$, $A$ and $B$ as input answers YES or NO to the question of whether they are separable:
1. IF $A$ and $B$ are not generically separable return "NO, $A$ and $B$ are not separable".
2. IF for some $W \in \mathfrak{L}_{M,A,B}$ the shadows $A_W$ and $B_W$ are not generically separable, return "NO, $A$ and $B$ are not separable".
3. ELSE return "YES, $A$ and $B$ are separable".

□

Dipartimento di Matematica, Università di Pisa, Via Buonarroti 2, 56127 Pisa, Italy
*E-mail address*: acquistf@gauss.dm.unipi.it

Departamento de Algebra, Facultad de Matemáticas, Universidad Complutense, 28040 Madrid, Spain
*E-mail address*: andradas@matss2.mat.ucm.es

Dipartimento di Matematica, Università di Pisa, Via Buonarroti 2, 56127 Pisa, Italy
*E-mail address*: broglia@gauss.dm.unipi.it